\documentclass[12pt]{article}
\usepackage{hyperref,amsfonts,amssymb,amsmath}
\usepackage{mathtools}

\usepackage{tikz-cd}
\usepackage{tikz-cd, amsmath, amssymb}

\usepackage{color}

\usepackage{soul}

\def\C{\mathbb{C}}

\def\Z{\mathbb{Z}}
\def\Q{\mathbb{Q}}

\def\R{\mathbb{R}}
\def\P{\mathbb{P}}
\def\A{\mathbb{A}}

\def\d{\partial}

\def\bq{ \begin{equation} }
\def\eq{ \end{equation} }
\def\ben{ \begin{eqnarray} }
\def\en{ \end{eqnarray} }

\def\frac#1#2{{#1\over #2}}

\def\on#1#2{\mathop{\vbox{\ialign{##\crcr\noalign{\kern2pt}
$\scriptstyle{#2}$\crcr\noalign{\kern2pt\nointerlineskip}
\kern-2pt$\hfil\displaystyle{#1}\hfil$\crcr}}}\limits}

\newtheorem{Constr}{Construction}

\begin{document}

\title{When the Fourier transform is one loop exact?}
\author{Maxim Kontsevich,  Alexander Odesskii}
   \date{}
\vspace{-20mm}
   \maketitle
\vspace{-7mm}
\begin{center}
IHES, 35 route de Chartres, Bures-sur-Yvette, F-91440,
France  \\[1ex]
and \\[1ex]
Brock University, 1812 Sir Isaac Brock Way, St. Catharines, ON, L2S 3A1, Canada\\[1ex]
e-mails: \\
\texttt{maxim@ihes.fr}\\
\texttt{aodesski@brocku.ca}
\end{center}

\medskip

\begin{abstract}

 We investigate the question: for which functions $f(x_1,...,x_n),~g(x_1,...,x_n)$ the asymptotic expansion of the integral $\int g(x_1,...,x_n) e^{\frac{f(x_1,...,x_n)+x_1y_1+...+x_ny_n}{\hbar}}dx_1...dx_n$ 
  consists only of the first term. We reveal a hidden projective invariance of the problem  which establishes its relation with geometry of projective hypersurfaces of the form $\{(1:x_1:...:x_n:f)\}$. We also construct various examples, in particular we prove that Kummer surface in $\P^3$ gives a solution to our problem.

\end{abstract}

\newpage

\tableofcontents

\newpage

 \section{Introduction}

 This paper is inspired by  beautiful results of \cite{efp} where a number of surprising identities for the Fourier transforms were discovered, for example    
 $$\int_{\R^2} sign(x_2)|x_2|^{-\frac{2}{3}}e^{i\frac{x_1^3}{x_2}}e^{i(x_1y_1+x_2y_2)}dx_1dx_2=\frac{2\pi i}{\sqrt{3}}\,sign(y_2)|y_2|^{-\frac{2}{3}}e^{\frac{i}{27}\frac{y_1^3}{y_2}}.$$
 The problem studied in \cite{efp} is the following. Let $F$ be a local field, $\psi$ a nontrivial unitary additive character of $F$, and $\chi_1,...,\chi_k$ multiplicative characters of $F$. A complex valued distribution of the form 
 $$\psi(Q(x_1,...,x_n))\prod_{j=1}^k\chi_j(P_j(x_1,...,x_n))$$
 is called elementary if $Q$ is a rational function and $P_j$ are polynomials. The problem is when the Fourier transform of such a distribution is also elementary. It was observed in \cite{efp} that one can set formally $\psi(x)=e^{\frac{ix}{\hbar}}$ and use the formal stationary phase method. This observation leads us to the following formal version of the problem where we no longer consider actual integrals and deal with the asymptotic expansion of oscillating integrals. Moreover, our functions $f,g$ (analogous of $Q,P_j$ from \cite{efp}) are no longer assumed to be rational and can be locally analytic or  just formal germs.
 
 Let $f(x_1,...,x_n)$ be a function in $n$ variables such that its Hessian is not identically zero:
 \begin{equation}\label{hess}
 \det\Bigg(\frac{\partial^2f}{\partial x_i\partial x_j}\Bigg)_{i,j}\ne 0.
 \end{equation}
 Let $\vec{x}=\vec{t}$ be a critical point of the function\footnote{Here and in the sequel we use vector notations like $\vec{x}=(x_1,...,x_n)$, $\vec{t}=(t_1,...,t_n)$ etc.} 
 \begin{equation}\label{fy}
 f(x_1,...,x_n)+x_1y_1+...+x_ny_n.
  \end{equation}
 It follows from (\ref{hess}) that the mapping $\vec{y}\mapsto\vec{t}$ has a non-degenerate Jacobian at generic point.  Indeed, equating to zero first derivatives of (\ref{fy}) 
 with respect to variables  $x_1,...,x_n$ we get
 \begin{equation}\label{yt}
 y_i=-\frac{\partial f(t_1,...,t_n)}{\partial t_i},~~~i=1,...,n
 \end{equation}
 and the Jacobian of this map is proportional to the Hessian (\ref{hess}) at $\vec{x}=\vec{t}$.
 Let $g(x_1,...,x_n)$ be another function in $x_1,...,x_n$. By the formal Fourier transform of the function $g(\vec{x})e^{\frac{f(\vec{x})}{\hbar}}$ we mean the perturbative 
 expansion of the formal integral  
 $$\int g(x_1,...,x_n) e^{\frac{f(x_1,...,x_n)+x_1y_1+...+x_ny_n}{\hbar}}dx_1...dx_n$$  
 given by the stationary phase method at the critical point $\vec{x}=\vec{t}$. Recall that  this 
 expansion has the form 
 \begin{equation}\label{per}
 \int g~ e^{\frac{f+x_1y_1+...+x_ny_n}{\hbar}}dx_1...dx_n=(2\pi\hbar)^{\frac{n}{2}}e^{\frac{\hat{f}}{\hbar}} \det\Bigg(-\frac{\partial^2f}{\partial x_i\partial x_j}\Bigg)_{i,j}^{-\frac{1}{2}}\Bigg\rvert_{\vec{x}=\vec{t}}~\sum_{k=0}^{\infty}A_k\hbar^k
 \end{equation}
 where coefficients $A_k$ of formal power series in $\hbar$ can be written as differential polynomials in $f(t_1,...,t_n),~g(t_1,...,t_n)$ with coefficients in $\Q$, divided by $\det\Big(\frac{\partial^2f}{\partial x_i\partial x_j}\Big)_{i,j}^{3k}\Big\rvert_{\vec{x}=\vec{t}}\,$, see Appendix. 
 
 We want to study the following question: 
 for which functions $f,g$ we have 
 \begin{equation}\label{sys0}
 A_k=0 \text{~~~~for all~~~~} k\ge 1.
 \end{equation}
 This condition can be written more explicitly in the form\footnote{Here and in the sequel such identities for integrals with 
 unspecified domain of integration mean that the r.h.s. is the perturbative expansion of the l.h.s. given by stationary phase method. Moreover, we often write the r.h.s. up to multiplication by a fourth root of unity $\pm 1,\pm i$.}
 \begin{equation}\label{eq0}
 \int g(\vec{x})e^{\frac{f(\vec{x})+\vec{x}\cdot\vec{y}}{\hbar}}d^n\vec{x}=(2\pi \hbar)^{\frac{n}{2}} \hat{g}(\vec{y})e^{\frac{\hat{f}(\vec{y})}{\hbar}}. 
 \end{equation}
Here 
\begin{equation}\label{fhat}
\hat{f}(\vec{y})=(f(\vec{x})+\vec{x}\cdot\vec{y})\Big\rvert_{\vec{x}=\vec{t}}
\end{equation}
is the Legendre transform of  function $f$, i.e.  the critical value of $f(\vec{x})+\vec{x}\cdot\vec{y}$, and 
\begin{equation}\label{ghat}
 \hat{g}(\vec{y})=g(\vec{x})\cdot \det\Bigg(-\frac{\partial^2f}{\partial x_i\partial x_j}\Bigg)_{i,j}^{-\frac{1}{2}}\Bigg\rvert_{\vec{x}=\vec{t}}.
\end{equation}

By a loose analogy with the terminology from physics literature we call the condition (\ref{eq0}) the 1-loop exactness of the formal Fourier transform.\footnote{The name is not totally precise as the Feynman graphs appearing in the expansion (\ref{per}) are not necessarily connected, see Appendix.}

 {\bf Definition 1.1.} The pair of functions $f,g$ is called admissible if the Hessian of $f$ is not identically zero (see (\ref{hess})), the function $g$ is not identically zero, and the identity (\ref{eq0}) holds.
 
 Notice that if we fix a function $f$, then the set of functions $g$ such that the pair $f,g$ is admissible or $g=0$, is a vector space, which we denote by $V_f$.
 
 {\bf Definition 1.2.} The rank of a function $f$ is the dimension of vector space of functions
 $$V_f=\{g~\vert ~f,g~ \text{is an admissible pair or}~ g=0\}.$$  
 The rank of an admissible pair $f,g$ is the rank of $f$.

 {\bf Definition 1.3.} A function $f$ is admissible if its rank is larger than zero. In other words, $f$ is admissible if there exists a (non-zero) function $g$ such that pair $f,g$ is admissible.

 Sometimes, if $f$ is fixed and clear from the context, we will call $g$ admissible if $f,g$ is an admissible pair.

\ 

 {\bf Remark 1.1.} First, in the case when $f$ is a concave function defined globally on $\R ^n$, and satisfying the condition 
 $$\lim_{|\vec{x}|\to\infty}\frac{f(\vec{x})}{|\vec{x}|}=-\infty~~~(\text{e.g.}~~~f(\vec{x})=-\frac{1}{2}\sum_{i=1}^nx_i^2~),$$
 and $g$ is also defined globally, 
 $$|g(\vec{x})|\leq C_1e^{C_2|\vec{x}|},$$
 the integral $\int_{\R^n}g(\vec{x})e^{\frac{f(\vec{x})+\vec{x} \cdot \vec{y}}{\hbar}}d\vec{x}$ is convergent for $\hbar>0$ and admits the asymptotic expansion as in 
 (\ref{per}).

 \

 In the general case where $f,g$ are germs of analytic functions, or even formal power series with coefficients in a field ${\bf k}$ of characteristic zero, by the formal integral we mean the expression (\ref{per}) where all terms $A_i$ are differential polynomials divided by an integer power of the Hessian  and therefore make sense, and the front factors are considered as formal symbols. 

 Also, in the usual Fourier transform one uses $\sqrt{-1}$ in the exponent. We omit it for our convenience, in order to simplify formulas. 

 In the sequel we will also replace $\det\Big(-\frac{\partial^2f}{\partial x_i\partial x_j}\Big)_{i,j}$ by $\det\Big(\frac{\partial^2f}{\partial x_i\partial x_j}\Big)_{i,j}$, which changes the result by a fourth root of unity. 

 \ 

 {\bf Remark 1.2.}  Here we explain an explicit procedure producing coefficients $A_i$ above, see Appendix for details. The calculation of the integral in the l.h.s. of (\ref{per}) near non-degenerate critical point can be  reduced (after some shift of variables $x_1,...,x_n$) to the following case.

 Let $F=F_2+F_3+...$ be a formal power series in $n$ variables $x_1,...,x_n$ with coefficients in a field ${\bf k}\supset\Q$ where $F_i$ are homogeneous polynomials of degree $i$, and $F_2$ is a non-degenerate quadratic form. We denote by $F^{\prime\prime}_2$ the corresponding symmetric matrix. Let $g=1+...\in {\bf k}[[x_1,...,x_n]]$ be 
 another power series (series $g$ starts with $1$ just for convenience). The formal integral 
 $$\frac{(\det F^{\prime\prime}_2)^{\frac{1}{2}}}{(2\pi \hbar)^{\frac{n}{2}}}\int g(\vec{x}) e^{\frac{F(\vec{x})}{\hbar}}dx_1...dx_n=1+...\in {\bf k}[[\hbar]]$$
 can be defined in the following way:

 First we rescale variables $x_i=\sqrt{\hbar}~\tilde{x}_i$. After that we have 
 $$\frac{F(\vec{x})}{\hbar}=F_2(\tilde{x}_1,...,\tilde{x}_n)+\hbar^{\frac{1}{2}}F_3(\tilde{x}_1,...,\tilde{x}_n)+\hbar^{\frac{2}{2}}F_4(\tilde{x}_1,...,\tilde{x}_n)+...$$
 and 
 $$g(\vec{x}) e^{\frac{F(\vec{x})}{\hbar}}=g(\hbar^{\frac{1}{2}}\vec{\tilde{x}}) e^{F_2(\vec{\tilde{x}})}e^{\hbar^{\frac{1}{2}}F_3(\vec{\tilde{x}})+\hbar^{\frac{2}{2}}F_4(\vec{\tilde{x}})+...}$$
 The r.h.s. of this formula can be expanded as power series in $\hbar^{\frac{1}{2}}$ and $\tilde{x}_1,...\tilde{x}_n$. In this way we reduced to computing the integrals of the form $\int he^{F_2}d\tilde{x}_1...d\tilde{x}_n$ where $h$ is a monomial in $\tilde{x}_1,...\tilde{x}_n$. This formal integral is defined to be zero if 
 degree of $h$ is odd, and as $\frac{1}{m!}\Delta^m_{F_2}(h)\int e^{F_2}d\tilde{x}_1...d\tilde{x}_n$ if degree of $h$ is $2m$. 
 Here $\Delta_{F_2}=-\frac{1}{2}\sum_{i,j} b_{ij}\partial_{\tilde{x}_i}\partial_{\tilde{x}_j}$ where $(b_{ij})=(F^{\prime\prime}_2)^{-1}$, and we formally declare $\int e^{F_2}d\tilde{x}_1...d\tilde{x}_n:=\frac{(2\pi \hbar)^{\frac{n}{2}}}{(\det F^{\prime\prime}_2)^{\frac{1}{2}}}.$

 Notice that the final expression $\sum_{i=0}^{\infty}A_i\hbar^i$ is a power series in $\hbar$ because integrals with monomials of odd degree are zero.

 \ 

 {\bf Remark 1.3.} If $f,g$ is an admissible pair, then the pair $\hat{f},\hat{g}$ of functions from the r.h.s. of (\ref{eq0}) is also admissible, because the Fourier transform is essentially an involution. Moreover, we have isomorphism of vector spaces $V_f\cong V_{\hat{f}}$ given by $g\mapsto\hat{g}$. In particular, $\dim V_f=\dim V_{\hat{f}}$ so $f$ and $\hat{f}$ have the same rank.

\ 

 {\bf Remark 1.4.} Let $\Sigma\subset\P^{n+1}$ be projective hypersurface locally defined by 
 \begin{equation}\label{locs}
 x_{n+1}=x_0f\Big(\frac{x_1}{x_0},...,\frac{x_n}{x_0}\Big).
 \end{equation}
 Then its projective dual hypersurface  $\widehat{\Sigma}\subset\P^{n+1}$
 is locally defined by 
  \begin{equation}\label{locd}
 y_0=y_{n+1}\hat{f}\Big(\frac{y_1}{y_{n+1}},...,\frac{y_n}{y_{n+1}}\Big)
 \end{equation}
 where $\hat{f}$ is the Legendre transform of $f$ given by (\ref{fhat}) and coordinates $(y_0,...,y_{n+1})$ are dual to $(x_0,...,x_{n+1})$. This projective duality plays an important role in this paper.

 \

 Let us describe the content of the paper. 

 In Section 2 we introduce the formalism of formal wave functions which allows to write expressions involving exponentials, integrals and delta functions in a purely algebraic context over an arbitrary field $\bf k$ of characteristic zero.

 In Section 3 we reformulate our problem in terms of  projectively dual hypersurfaces defined by (\ref{locs}) and (\ref{locd}). We show that admissibility condition does not depend on the choice of projective coordinates. In this way we deduce a projective invariance of the original problem: if $f$ is admissible and we make an arbitrary linear change of variables 
 $$x_i=\sum_{j=0}^{n+1}a_{i,j}\tilde{x}_j,~~~i=0,...,n+1$$ 
 in (\ref{locs}), solve with respect  to  $\tilde{x}_{n+1}$ writing the result in the form
 \begin{equation}\label{act}
  \tilde{x}_{n+1}=\tilde{x}_0\tilde{f}\Big(\frac{\tilde{x}_1}{\tilde{x}_0},...,\frac{\tilde{x}_n}{\tilde{x}_0}\Big),
  \end{equation}
 then $\tilde{f}$ is also admissible for arbitrary non-degenerate matrix $(a_{i,j})\in GL(n+2)$.

 In Section 4 we explain how to reduce our problem in the case of algebraic hypersurfaces to a finite system of differential equations for $f,g$. Notice that in general our system of differential equations (\ref{sys0}) is infinite. In the special case when $g$ is a solution of a regular holonomic system of differential equations, our question about admissible pairs essentially reduces to the following interesting problem concerning algebraic holonomic $D$-modules: 
 
 {\it find regular holonomic $D$-modules $M$ on the affine space $\A^{n+2}$ which are monodromic (i.e. the action of the Euler operator $\sum_{i=0}^{n+1}x_i\partial_{x_i}$ is locally finite) and such that the singular support of $M$ does not contain conormal bundles to $\A^{n+2}$ and $\{0\}\subset\A^{n+2}$.}

 \

 In Section 5 we study numerous examples. In particular, we prove that Kummer quartic surface in $\P^3$ and Segre cubic in $\P^4$ are both admissible. We also construct a huge family of admissible functions for arbitrary $n\geq 2$. These functions are defined parametrically as 
 $$f=\frac{\phi_{n+1}}{\phi_0},~~~x_i=\frac{\phi_i}{\phi_0},~~~i=1,...,n,$$
 $$\phi_i=\frac{1}{2}\sum_{j,k=0}^na_{i,j,k}u_ju_k,~~~i=0,...,n+1$$
  where $u_0,...,u_n$ are coordinates on the corresponding hypersurface, $a_{i,j,k}\in {\bf k}$ are constants and $a_{i,k,j}=a_{i,j,k}.$ We prove that any non-degenerate function of this family (i.e. with the generically non-vanishing Hessian) is admissible and $\dim V_f\geq n+1$ for $n>2$. We also construct examples of such functions for arbitrary $n\geq 2$ with the rank equal to $(n+1)!$.

 In Section 6 we present the classification of admissible pairs of functions in one variable. 
 
 In Section 7 we show some classification results of admissible pairs of functions in two variables. 

 In Section 8 we outline a potential application of our studies to the construction of generalized Dirichlet series based on the Poisson summation formula and the Mellin transform.

 In Section 9 we formulate several conjectures and open questions. Notice that some other conjectures and open questions are discussed in the main part of the paper. 

 In Appendix we recall how to write explicitly the infinite system of differential equations (\ref{sys0}) for $f,g$ based on the stationary phase method and the Feynman diagrams technique.

 \section{Formal wave functions}
 
 Here we introduce a rigorous language of wave functions which allows us to use expressions like exponential functions depending on small parameter $\hbar$ and written as  $\exp(f(x_1,\dots,x_n)/\hbar)$, or delta functions $\delta(f(x_1,\dots, x_n))$ etc,  in a purely  algebraic situation, when neither coordinates nor functions take real or complex values. The whole calculus makes sense over any field $\bf k$ of characteristic zero.
 
 Let $M=\A^{2n}$ be an affine space of dimension $2n$ over $\bf k$, endowed with a translationally-invariant symplectic structure. In other words, we have global coordinates\footnote{In the sequel we will often use notation $y_1=x_{n+1},...,y_n=x_{2n}$.} $(x_1,\dots,x_{2n})$ on $M$ defined up to affine-symplectic  transformations
 $$ x_i\mapsto \sum_{j=1}^{2n} a_{ij} x_j +b_j $$
 where $(a_{ij})_{1\le i,j,\le 2n}$ is an invertible matrix preserving the standard symplectic 2-form
 $$\omega=\sum_{i=1}^n dx_i\wedge d x_{i+n}\,.$$
 Denote by $(\gamma_{ij})_{1\le i,j\le 2n}$ the tensor for the inverse bi-vector field:
 $$\gamma_{ij}:=\begin{cases}1& \mbox{ if } j=i+n,\quad 1\le i\le n,\\
 -1& \mbox{ if } i=j+n,\quad 1\le j\le n,\\
 0&\mbox{ otherwise}\end{cases}   $$
 
 Let $\hbar$ be a formal variable. We define the canonical Moyal star-product on the vector space
 $${\bf k}[x_1,\dots,x_{2n},\hbar] $$ by the formula
 $$f\star g:= \sum_{k=0}^\infty \frac{\hbar^k}{k!} \left[\left(\sum_{1\le i,j\le 2n} \frac{\gamma_{ij}}{2} \frac{\d}{\d x_i^{(1)}}\frac{\d}{\d x_j^{(2)}}\right)^k(f((x_i^{(1)})_{1\le i\le 2n} \cdot g((x_i^{(2)})_{1\le i\le 2n}) \right]\Bigg\vert_{x_i^{(1)}=x_i^{(2)}=x_i,\,\,\,\, i=1,\dots ,2n}$$
 where in the square brackets we consider functions on $M\times M$ endowed coordinates
 $$ (x_1^{(1)}, \dots,  x_{2n}^{(1)}),~ (x_1^{(2)}, \dots,  x_{2n}^{(2)})$$
 (i.e. two copies of the original coordinates on $M$).
 One can rewrite the above formula as
 $$  f\star g=\left[\exp\left(\frac{1}{\hbar} \sum_{1\le i,j\le 2n}\frac{1}{2}\gamma_{ij} \frac{\d}{\d x_i^{(1)}}\boxtimes \frac{\d}{\d x_j^{(2}}\right)(f\boxtimes g)\right]\Bigg\vert_{\text{diagonal}} \,.$$

 Algebra $({\bf k}[x_1,\dots,x_{2n},\hbar],
\star) $ over $\bf k$ is an {\it associative unital} algebra generated by $(2n+1)$  elements $(\hat{x}_i)_{i=1,\dots,2n}, \hbar$
  where $\hat{x}_i$ correspond to $x_i\in {\bf k}[x_1,\dots,x_{2n},\hbar]$, satisfying the relations
  $$[\hat{x}_i,\hat{x}_j]=\gamma_{ij}\hbar, \quad [\hat{x}_i,
  \hbar]=0$$
  and can be identified with the algebra\footnote{We will see later that it is more natural to consider differential operators acting on {\it half-densities} instead of functions.} of polynomial $\hbar$-differential operators in $n$ variables $x_1,\dots,x_n$ by
  $$\hat{x}_i\mapsto x_i, \quad \hat{x}_{n+i}\mapsto \hbar \frac{\partial}{\partial {x_i}},\qquad \forall i=1,\dots,n\,. $$

 The Moyal product is covariant with respect to the action of the group of affine-symplectic automorphisms of $M$. It extends to the ${\bf k}[[\hbar]]$-linear product on $A[[\hbar]]$ where $A$ is the algebra of functions on a Zariski open subset of $M$, or analytic functions in a Stein open domain in $M$ if ${\bf k}=\C$, or $C^\infty$-functions in an open domain if ${\bf k}=\R$. Also, for general field $\bf k$ of characteristic zero, one can take $A$ to be the algebra of formal power series at a given point $m\in M$.
 
 \
 
 Let $L$ be a Lagrangian submanifold of $M$ in a broad sense, i.e. an algebraic subvariety, or analytic/smooth in the case ${\bf k}=\C$ or ${\bf k}=\R$, or a germ of such subvariety, or even a {\it formal germ} at some point of $M$.
 
  We assume that a spin-structure on $L$ is given, which means that we are given a line bundle on $L$ whose tensor square is identified with the canonical bundle $K_L$, i.e. the bundle $\wedge^n T^*_L$ of top-degree forms on $L$. The chosen square root bundle we denote by $K_L^{\otimes 1/2}$. In the case of a germ, or a formal germ at point $m\in M$, it is sufficient to choose a square root of the fiber $(K_L)_{Im}$ at the base point $m$.

  Our main goal is the following
  \begin{Constr}\!\!\!.
   With a given pair $(L,K_l^{\otimes 1/2})$ we associate a module (or more, precisely, a sheaf of modules)  ${\mathcal WF}_{L,K_L^{\otimes 1/2}}$ over the quantum algebra $A[[\hbar]]$ where $A$ is the algebra of functions in the formal completion of $M$ at $L$.
  \end{Constr}
  
  Elements of the vector space  ${\mathcal WF}_{L,K_L^{\otimes 1/2}}$ are called  {\it formal wave functions supported on $L$}.
  
  There are several approaches to this constructions. The one presented below is an explicit one in affine symplectic coordinates, but it is based on a nontrivial consistency check\footnote{There is another approach (not described in this paper) based on Gelfand-Kazhdan type formal geometry.}.

  First, we give the definition of  the module  ${\mathcal WF}_{L,K_L^{\otimes 1/2}}$ in local coordinates. Let $(x_1,\dots,x_{2n})$ be  global affine symplectic coordinates on $M$ such that the symplectic form is the standard one, and such that locally near a point of $L$ the projection $\pi$  to the coordinate space $\mathbb A^n$ by first $n$ coordinates $(x_1,\dots,x_n)$  has  a {\it non-zero} Jacobian. Hence, locally $L$ is a graph of a closed 1-form on $\mathbb A^n$:
   $$x_{n+i}=\alpha_i(x_1,\dots,x_n)\quad\forall i=1,\dots,n,\qquad d\alpha=0,\quad \alpha:=\sum_{i=1}^n \alpha_i d x_i \,.$$
Let us also choose a generator of $K_L^{\otimes 1/2}$ whose tensor square is $\pi^*(dx_1\wedge \dots\wedge dx_n)$. We will denote in short this generator by $(dx_1\wedge \dots\wedge dx_n)^{1/2}$, it is well-defined up to a sign.

After we make the choices from above (i.e. affine symplectic coordinate system $(x_1,\dots,x_{2n})$ and the generator $(dx_1\wedge \dots\wedge dx_n)^{1/2}$), we declare 
${\mathcal WF}_{L,K_L^{\otimes 1/2}}$ to be the ${\bf k}[[\hbar]]$-module
   $$\mathcal O(L)[[\hbar]]\simeq \mathcal O[\pi(L)][[\hbar]]$$
   where $\mathcal O(L)$ is the algebra of functions on $L$, identified with functions on open domain (or a formal germ) $\pi(L)\subset \mathbb A^n$, with the action of generators $\hat{x}_i,\quad i=1,\dots ,2n$ of the quantum algebra given by
   $$\hat{x}_i\mapsto x_i,\quad \hat{x}_{n+i}\mapsto \hbar \frac{\partial}{\partial {x_i}}+\alpha_i,\qquad \forall i=1,\dots,n\,. 
   $$
   It is a nontrivial fact that this action extends by continuity from polynomial elements to the functions on the affine space $\mathbb A^{2n}\simeq M$ defined in the formal completion of $L\subset M$.

 Notationally, it is convenient to choose locally a primitive of closed 1-form $\alpha$, i.e. a function $f(x_1,\dots,x)$ defined on the domain $\pi(L)\subset \mathbb A^n$, such that 
 $$\alpha=df\iff \alpha_i=\frac{\partial f}{\partial x_i} \quad \forall i=1,\dots,n\,.$$
 Then the element of ${\mathcal WF}_{L,K_L^{\otimes 1/2}}$ corresponding to 
 a series  
 $$g=g(x_1,\dots,x_n;\hbar)=\sum_{k\ge 0}g_k \hbar^k\in \mathcal O[\pi(L)][[\hbar]]$$
 we denote by the formal product 
  $$\psi= e^{\frac{f(x_1,\dots,x_n)}{\hbar}}g(x_1,\dots,x_n;\hbar)\cdot  (dx_1\wedge \dots\wedge dx_n)^{1/2} \,.$$
  The choice of the primitive is irrelevant: if we shift it by a constant $f\to f+const$ then we formally multiply the expression above by $e^{\frac{const}{\hbar}}$ without affecting  series $g$.
  
  \
  
  Now we want to study the dependence of the description of ${\mathcal WF}_{L,K_L^{\otimes 1/2}}$ under the change of choices made.
  First, for a given affine symplectic coordinate system $(x_1,\dots,x_{2n})$, if we change the square root $(dx_1\wedge \dots\wedge dx_n)^{1/2}$ by sign, then the isomorphism
   $${\mathcal WF}_{L,K_L^{\otimes 1/2}}\stackrel{\sim}{\longrightarrow} O[\pi(L)][[\hbar]]$$
  will be also changed by sign.
  
  Next, let us change the affine symplectic coordinate system $(x_1,\dots,x_{2n})$ preserving first $n$ coordinates $x_1,\dots,x_n$:
  \begin{equation}  x_i\to x_i, \quad x_{n+i}\to x_{n+i}+\sum_{j=1}^n b_{ij} x_j+c_i,\qquad i=1,\dots,n,\quad b_{ij}=b_{ji}\in \mathbf k,\quad c_i \in \mathbf k\,. \label{change1}\end{equation}
  Then we multiply the corresponding formal wave function $\psi$ by
  $$\psi \to e^{\frac{1}{\hbar}\left(\frac{1}{2}\sum_{i,j} b_{ij} x_i x_j +\sum_i c_i x_i\right)}\psi $$
  If we apply an affine symplectic transformation associated with an invertible $(n\times n)$-matrix $a=(a_{ij})_{1\le i,j\le n}$:
  \begin{equation}x_i\to x_i'=\sum_{j=1}^n a_{ij} x_j, \quad x_{n+i}\to x_{n+i}'=\sum_{j=1}^n (a^{-1})_{ji} x_{n+j} \qquad \forall i=1,\dots,n  \label{change2}\end{equation}
  then we change the wave function by
    $$ e^{\frac{1}{\hbar} f} g\to   e^{\frac{1}{\hbar} f'} g',\qquad f'(\vec x):=f(\vec x\,'(\vec x)),\quad g'(\vec x;\hbar):=g(\vec x\,'(\vec x);\hbar)\cdot \det(a)^{-1/2}\,.$$
    
    We can now describe the dependence of the description of ${\mathcal WF}_{L,K_L^{\otimes 1/2}}$ under the (almost) general affine symplectic transformation. Namely, let us assume that $(x_1,\dots,x_{2n})$ and $(x_1',\dots,x_{2n}')$ are two affine symplectic coordinate systems on $M$ such that both projections from $L$ to the affine space $\mathbb A^n$ given by $(x_1,\dots,x_n)$ and  $(x_1',\dots,x_n')$ are open embeddings.
     Let us make an assumption (which is an open condition) that
      $(x_1,\dots,x_n,x_1',\dots,x_n')$ form a system of coordinates on $M$. By applying the above modifications \eqref{change1},\eqref{change2}, we may assume that
      $$x_i'=x_{n+i},\quad x_{n+i}'=-x_i,\qquad \forall i=1,\dots,n $$
      Then we declare that the corresponding formal wave functions undergo the {\it formal Fourier transform}:
      $$ g(\vec x;\hbar) e^{\frac{f(\vec x)}{\hbar}}\to \tilde{g}(\vec y;\hbar) e^{\frac{\tilde{f}(\vec y)}{\hbar}} :=\frac{1}{(2\pi \hbar)^{n/2}}\int g(\vec x;\hbar) e^{\frac{f(\vec x)+\vec x\cdot \vec y}{\hbar}} d^n \vec x $$
      where the integral in the r.h.s. is understood as the asymptotic expansion, calculated via the stationary phase method. In particular, the exponent $\tilde 
      f$ is the Legendre transform of $f$:
      $$\hat{f}(\vec y)=\text{Critical value of }  f(\vec x)+\vec x\cdot \vec y
      \,.
      $$
      The consistency check mentioned  before, says that for {\it three} 
        affine symplectic coordinate systems on $M$: $(x_1,\dots,x_{2n})$, $(x_1',\dots,x_{2n}')$ and $(x_1'',\dots,x_{2n}'')$, the passage from the first to the second, and then from the second to the third, coincides with the passage from the first to the  third. Here we will give the sketch of the proof which is not purely algebraic and is  based partially on analysis.

        After making the assumption that the triple of coordinate systems under consideration is sufficiently generic,  the consistency question can be reduced to the following equality
        \begin{equation}\int \Bigg(e^{-{\vec y\cdot \vec y\over 2\hbar}}\Big(\int g(\vec x) e^{f(\vec x)+\vec x\cdot\vec y \over \hbar}d^n\vec x\Big)\Bigg)e^{\vec y\cdot \vec z\over \hbar}d^n\vec y=(2\pi\hbar)^{n/2}e^{{\vec z\cdot\vec z\over 2\hbar}}\int e^{{1\over \hbar}(f(\vec x)+{\vec x\cdot \vec x\over 2}+\vec x\cdot \vec z)}d^n\vec x\label{consistency}\end{equation}
        where $f(\vec x)$ is a formal power series starting with quadratic terms
        \begin{equation}f(\vec x)=-{1\over 2}\sum_{ij} b_{ij}x_i x_j+\dots \label{quadratic_part} \end{equation}
        where symmetric matrix $(b_{ij})_{1\le i,j\le n}$ has no eigenvalues equal to $0$ or $1$.
        
        In the case $\mathbf k=\mathbb R$, function $f$ being a global strictly concave smooth real-valued function on $\R^n$ such that $f(\vec x)\le -{c\over 2}\vec x\cdot \vec x$ for some $c>1$ and any smooth function $g$ of at most exponential growth at infinity, all the integrals in \eqref{consistency} are absolutely convergent. The equality follows because the l.h.s can be rewritten as
        \begin{equation*}\int\int g(\vec x) e^{{1\over \hbar}(f(\vec x)+\vec x\cdot \vec y+\vec y\cdot \vec z-{\vec y\cdot \vec y \over 2})}d^n\vec x \,d^n\vec y\end{equation*}
        and then identified with r.h.s. using the rewriting
        $$\vec x\cdot \vec y+\vec y\cdot \vec z-{\vec y\cdot \vec y \over 2}=-{\vec w\cdot\vec w\over 2}+{\vec x\cdot\vec x\over 2}+{\vec z\cdot \vec z\over 2}+\vec x\cdot\vec z,\quad \vec w:=\vec y-\vec x-\vec z \,.$$
        Each term in the $\hbar$-expansion of the equality \eqref{consistency} is a polynomial identity with rational coefficients involving {\it finitely many} Taylor coefficients of $f,g$. The fact that it holds for real $C^\infty$ examples as above giving  Zariski dense subsets of possible Taylor coefficients of finite jets of $f,g$ at $0$ implies that \eqref{consistency} holds for arbitrary formal series $f,g$ with coefficients in any field  $\mathbf k\supset \mathbb Q$.  This concludes the proof of the consistency check.

         \

        {\bf Remark 2.1.} One can further generalize equality \eqref{consistency}. Namely, each term in $\hbar$-expansion is an equality of certain finite sums of numbers obtained by contraction of upper and lower indices for certain symmetric tensors in $n$-dimensional space. The tensors under consideration are Taylor coefficients of series $f,g$ and the inverses to symmetric matrices $B$ and  $B-{\bf 1}_n$
         where $B=(b_{ij})_{1\le i,j\le n}$ is the (negative) Hessian of $f$ at $0$ (see \eqref{quadratic_part}) and ${\bf 1}_n$ is the identity matrix (could be replaced by any non-degenerate quadratic form).  The fact that the equality holds in {\it any} positive dimension $n\ge 0$ implies by the Weyl's  fundamental theorem in invariant theory that it  holds by purely formal reasons, as the cancellation of linear combination of oriented graphs controlling the contraction of indices. Therefore, the equality \eqref{consistency} makes sense and holds in arbitrary $\mathbb Q$-linear rigid symmetric monoidal category, like, e.g. finite-dimensional super vector spaces. Hence the construction $$(L,K_L^{\otimes 1/2})\leadsto {\mathcal W }F_{L,K_L^{\otimes 1/2}}$$ can be extended to the case of super manifolds.

         \

        In what follows we will not specify the choice of $K_L^{\otimes 1/2}$ and hence omit it from the notation.

      \

      {\bf Remark 2.2.} Let $f,g$ be an admissible pair. The formal product 
      $$\psi=g(x_1,...,x_n)e^{\frac{1}{\hbar}f(x_1,...,x_n)}(dx_1\wedge ...\wedge dx_n)^{\frac{1}{2}}$$
      is an element of ${\mathcal WF}_{L}$ where $L\subset\A^{2n}$ is the graph of $df$, a germ of Lagrangian submanifold in $\A^{2n}=T^*\A^n$. Simultaneously, $L\subset\A^{2n}=T^*(\A^n)^*$ is the graph of $d\hat{f}$ where $\hat{f}(y_1,...,y_n)$ is the Legendre transform of $f$, see (\ref{fhat}). 

      The admissible pair $\hat{f},\hat{g}$ obtained by the Fourier transform from $f,g$ gives {\it the same} element $\hat{\psi}=\psi\in{\mathcal WF}_{L}$.
      Our initial question about finding admissible pairs can be reformulated as the question about finding elements $\psi\in{\mathcal WF}_{L}$ such that in two descriptions of ${\mathcal WF}_{L}$ corresponding to the projections either to coordinates $x_1,...,x_n$ or to $x_{n+1},...,x_{2n}$, the functions $g(x_1,...,x_n,\hbar)$ (and $\hat{g}(x_{n+1},...,x_{2n},\hbar)$) associated with $\psi=\hat{\psi}$ do not depend on $\hbar$.

         \ 
      
      Finally, we explain how to interpret in terms of formal wave functions expressions involving delta functions.
       Let us assume that we have a (germ) of $k$-dimensional submanifold in a  $n$-dimensional affine space. After making an affine change of coordinates, we may assume that the submanifold under the consideration is given the graph of a map form an open domain in $\mathbb A^k$ to $\mathbb A^{n-k}$, i.e. given by
        \begin{align*}x_{k+1}&=\phi_1(x_1,\dots,x_k)\\
          x_{k+2}&=
          \phi_2(x_1,\dots,x_k)\\
          \dots\\
         x_{n}&=\phi_{n-k}(x_1,\dots,x_k) \,.\end{align*}
         
         Assume that we are also given two functions $f(x_1,\dots,x_k)$ and $g(x_1,\dots,x_k)$. We would like to make sense of the following expression:
         $$ g(x_1,\dots,x_k) e^{\frac{1}{\hbar} f(x_1,\dots,x_k)}\prod_{i=1}^{n-k} \delta(x_{k+i}- \phi_i
         (x_1,\dots,x_k) ) \,.$$
         which is an element of ${\mathcal WF}_{L}$ where $L\subset M\simeq\A^{2n}=T^*\A^n$ is the conormal bundle to the $k$-dimensional submanifold in $\A ^n$ defined above.
         This can be achieved by making the Fourier transform in variables
          $(x_{k+1},\dots,x_n)$:
          \begin{multline*}
              \int g(x_1,\dots,x_k) e^{\frac{1}{\hbar} (f(x_1,\dots,x_k)+\sum_{i=1}^{n-k} x_{k+i} y_i)}\prod_{i=1}^{n-k} \delta(x_{k+i}- \phi_i
         (x_1,\dots,x_k) ) \prod_{i=1}^{n-k} d x_{k+i}=\\
         =g(x_1,\dots,x_k) e^{\frac{1}{\hbar} (f(x_1,\dots,x_k)+\sum_{i=1}^{n-k} \phi_i(x_1,\dots,x_k) y_i)}\,.
          \end{multline*} 
 
  So, we see that new exponent which is function in variables  $x_1,\dots,x_k;y_1,\dots,y_{n-k}$
  $$ f(x_1,\dots,x_k)+\sum_{i=1}^{n-k} \phi_i(x_1,\dots,x_k) y_i   $$
  which happen to be  a linear function in $(n-k)$ variables $y_1,\dots,y_{n-k}$.

  In our formalism we have
  \begin{equation}\label{del}
  \delta(x)=\frac{1}{2\pi  \hbar}\int e^{\frac{xy}{\hbar}}dy.
  \end{equation}
  Recall that all our identities hold up to fourth root of unity (see Introduction). In the case of the actual Dirac distribution $\delta(x)$ on $\R$, the exact formula is 
  $$
  \delta(x)=\frac{1}{2\pi  \hbar}\int_{-\infty}^{\infty} e^{\frac{xy}{i\hbar}}dy,~~~\hbar>0.
  $$

   \ 

  {\bf Remark 2.3.} Algebra ${\bf k}[x_1,...,x_{2n},\hbar]$ with Moyal star product (as well as its completions associated with open domain in $\A^{2n}$ or formal completions) has a canonical derivation over $\bf k$ given by 
  $$\tau(x_i)=\frac{1}{2}x_i,~~~\tau(\hbar)=\hbar.$$
  If $L\subset\A^{2n}$ is conical, then $\tau$ admits a natural extension $\tau_L$ to ${\mathcal WF}_{L}$. Locally, if $x_1,...,x_n$ are coordinates, then 
  $$L=\text{graph}~dF_L(x_1,...,x_n)$$
  where $F_L$ is homogeneous of degree 2. Let
  $$\psi=G_L(x_1,...,x_n,\hbar)e^{\frac{F_L(x_1,...,x_n)}{\hbar}}(dx_1...dx_n)^{\frac{1}{2}}\in{\mathcal WF}_{L}.$$
  We define 
  $$\tau_L(\psi)=(\hbar\partial_{\hbar}G_L(x_1,...,x_n,\hbar)+\frac{1}{2}\sum_{i=1}^nx_i\partial_{x_i}G_L(x_1,...,x_n,\hbar))e^{\frac{F_L(x_1,...,x_n)}{\hbar}}(dx_1...dx_n)^{\frac{1}{2}}.$$

  \section{Reformulation of the problem in terms of conical germs}

  \subsection{Admissible pairs and projectively dual hypersurfaces}

  {\bf Definition 3.1.1.} A germ of smooth hypersurface $\Sigma$ in projective space $\P (V)$ is called non-degenerate if the Gauss map $\Sigma\to\P(V^*)$ given by $x\mapsto T_x\Sigma$ is an immersion.

  Projective duality identifies germs of smooth non-degenerate hypersurfaces in $\P(V)$ and $\P(V^*)$. 

  {\bf Theorem 3.1.1.} There is one to one correspondence between admissible pairs of germs of functions $f,g$ in $n$ variables and germs of distributions in $n+2$ variables which are smooth densities on a conical germ of a hypersurface in $\A ^{n+2}$, independent of $\hbar$, satisfying certain genericity constraints explained below, homogeneous of degree $-\frac{n+2}{2}$ and such that their Fourier transform is (up to the formal factor $(2\pi\hbar)^{\frac{n+2}{2}}$) a distribution with the same properties on the dual space. This correspondence is given by 
  \begin{equation}\label{G}
  (f,g)~~~\leftrightsquigarrow ~~~G(x_0,...,x_{n+1})=\delta\Big(\frac{x_{n+1}}{x_0}-f\Big(\frac{x_1}{x_0},...,\frac{x_n}{x_0}\Big)\Big)~g\Big(\frac{x_1}{x_0},...,\frac{x_n}{x_0}\Big)~x_0^{-\frac{n+2}{2}}
  \end{equation}
  and we have
  \begin{equation}\label{FG}
  \int G(x_0,...,x_{n+1})e^{\frac{1}{\hbar}(x_0y_0+...+x_{n+1}y_{n+1})}dx_0...dx_{n+1}=(2\pi\hbar)^{\frac{n+2}{2}}\hat{G}(y_0,...,y_{n+1})
  \end{equation}
  where $\hat{G}$ is given by (\ref{Ghat}) for some germs of functions $\hat{f},\hat{g}$.

  The genericity constraints\footnote{For any non-degenerate hypersurface in $\P(\A^{n+2})$ the corresponding germ of cones satisfies above constraints at generic point.} for a smooth conical germ $C$ at point $p\in\A^{n+2}\backslash  0$ are the following:

  {\bf 1.} The projectivization of $C$ is non-degenerate,

  {\bf 2.} Coordinate $x_0$ of point $p$ is non-zero,

  {\bf 3.} The tangent space $T_{x_0}C$ does not contain vector $(0,...,0,1)$.

  {\bf Proof} Let $f,g$ be admissible. Notice that the formula (\ref{G}) for $G$ describes a germ of a general smooth distribution in $n+2$ variables supported on a conical germ of a hypersurface in $\A ^{n+2}$, and homogeneous of degree $-\frac{n+2}{2}$. Let us check that Fourier transform of $G$ can be represented by a similar formula multiplied by $(2\pi\hbar)^{\frac{n+2}{2}}$. Indeed, we have 
  $$\int G(x_0,...,x_{n+1})e^{\frac{1}{\hbar}(x_0y_0+...+x_{n+1}y_{n+1})}dx_0...dx_{n+1}\stackrel{(1)}{=}$$
  $$\int \delta\Big(x_{n+1}-x_0f\Big(\frac{x_1}{x_0},...,\frac{x_n}{x_0}\Big)\Big)~g\Big(\frac{x_1}{x_0},...,\frac{x_n}{x_0}\Big)~x_0^{-\frac{n}{2}}e^{\frac{1}{\hbar}(x_0y_0+...+x_{n+1}y_{n+1})}dx_0...dx_{n+1}\stackrel{(2)}{=}$$
  $$\int g\Big(\frac{x_1}{x_0},...,\frac{x_n}{x_0}\Big)~x_0^{-\frac{n}{2}}e^{\frac{1}{\hbar}(x_0y_{n+1}f\big(\frac{x_1}{x_0},...,\frac{x_n}{x_0}\big)+x_0y_0+...+x_{n}y_{n})}dx_0...dx_n\stackrel{(3)}{=}$$
  $$\int g\big(x_1,...,x_n\big)~x_0^{\frac{n}{2}}e^{\frac{x_0y_0}{\hbar}+\frac{x_0y_{n+1}}{\hbar}\big(f(x_1,...,x_n)+x_1\frac{y_1}{y_{n+1}}+...+x_{n}\frac{y_{n}}{y_{n+1}}\big)}dx_0...dx_n\stackrel{(4)}{=}$$
  $$\int x_0^{\frac{n}{2}} e^{\frac{x_0y_0}{\hbar}}~\Bigg(\frac{2\pi\hbar}{y_{n+1}x_0}\Bigg)^{\frac{n}{2}}\hat{g}\Big(\frac{y_1}{y_{n+1}},...,\frac{y_n}{y_{n+1}}\Big)e^{\frac{y_{n+1}x_0}{\hbar}\hat{f}\Big(\frac{y_1}{y_{n+1}},...,\frac{y_n}{y_{n+1}}\Big)}dx_0\stackrel{(5)}{=}$$
  $$(2\pi\hbar)^{\frac{n+2}{2}}\delta\Big(y_0+y_{n+1}\hat{f}\Big(\frac{y_1}{y_{n+1}},...,\frac{y_n}{y_{n+1}}\Big)\Big)\hat{g}\Big(\frac{y_1}{y_{n+1}},...,\frac{y_n}{y_{n+1}}\Big)y_{n+1}^{-\frac{n}{2}}\stackrel{(6)}{=}$$
  $$(2\pi\hbar)^{\frac{n+2}{2}}\delta\Big(\frac{y_0}{y_{n+1}}+\hat{f}\Big(\frac{y_1}{y_{n+1}},...,\frac{y_n}{y_{n+1}}\Big)\Big)\hat{g}\Big(\frac{y_1}{y_{n+1}},...,\frac{y_n}{y_{n+1}}\Big)y_{n+1}^{-\frac{n+2}{2}}=$$
  $$(2\pi\hbar)^{\frac{n+2}{2}} \hat{G}(y_0,...,y_{n+1})$$
  where
  \begin{equation}\label{Ghat}
  \hat{G}(y_0,...,y_{n+1})=\delta\Big(\frac{y_0}{y_{n+1}}+\hat{f}\Big(\frac{y_1}{y_{n+1}},...,\frac{y_n}{y_{n+1}}\Big)\Big)\hat{g}\Big(\frac{y_1}{y_{n+1}},...,\frac{y_n}{y_{n+1}}\Big)y_{n+1}^{-\frac{n+2}{2}}
  \end{equation}
  and this formula is indeed similar to (\ref{G}).
  
  Here is an explanation of all steps in this calculation:

  (1) We use the homogeneous property of delta function.

  (2) We integrate by $x_{n+1}$ removing delta function.

  (3) We make a change of variables $x_i\mapsto x_0x_i,~i=1,...,n$.

  (4) We integrate by $x_1,...,x_n$ using our assumption that $f,g$ is admissible, see (\ref{eq0}).

  (5) We integrate by $x_0$ using an integral representation of delta function, see (\ref{del}).

  (6)  We use the homogeneous property of delta function again.

   \
   
  Conversely, let us assume that the Fourier transform of $G$ is supported on a conical germ and also is homogeneous of degree $-\frac{n+2}{2}$, i.e.  
  $$\int \delta\Big(x_{n+1}-x_0f\Big(\frac{x_1}{x_0},...,\frac{x_n}{x_0}\Big)\Big)~g\Big(\frac{x_1}{x_0},...,\frac{x_n}{x_0}\Big)~x_0^{-\frac{n}{2}}e^{\frac{1}{\hbar}(x_0y_0+...+x_{n+1}y_{n+1})}dx_0...dx_{n+1}=$$
  $$(2\pi\hbar)^{\frac{n+2}{2}}\delta\Big(y_0+y_{n+1}\hat{f}\Big(\frac{y_1}{y_{n+1}},...,\frac{y_n}{y_{n+1}}\Big)\Big)\hat{g}\Big(\frac{y_1}{y_{n+1}},...,\frac{y_n}{y_{n+1}}\Big)y_{n+1}^{-\frac{n}{2}}$$
  where $f,g,\hat{f},\hat{g}$ are some germs of functions in $n$ variables independent of $\hbar$. Integrating the l.h.s. by $x_{n+1}$ and making change of variables  $x_i\mapsto x_0x_i$, $~~y_i\mapsto y_{n+1}y_i,~i=1,...,n$ we obtain:
  $$\int g\big(x_1,...,x_n\big)~x_0^{\frac{n}{2}}e^{\frac{x_0y_0}{\hbar}+\frac{x_0y_{n+1}}{\hbar}(f(x_1,...,x_n)+x_1y_1+...+x_{n}y_{n})}dx_0...dx_n=$$
  $$(2\pi\hbar)^{\frac{n+2}{2}}\delta\Big(y_0+y_{n+1}\hat{f}\Big(y_1,...,y_n\Big)\Big)\hat{g}\Big(y_1,...,y_n\Big)y_{n+1}^{-\frac{n}{2}}.$$
  Finally, we multiply this equation by $\frac{1}{2\pi\hbar}e^{-\frac{y_0}{\hbar y_{n+1}}}$, integrate by $y_0$, and after subsequent integration of the l.h.s. by $x_0$ we obtain (\ref{eq0}). $\square$

  \ 

  {\bf Remark 3.1.1.} One can check that if the densities $G$ and $\hat{G}$ are related by the Fourier transform (\ref{FG}), then their supports are projectively dual. See also Remark 1.4.

  \subsection{Projective invariance of the problem}

  Notice that the conditions on distribution $G$ from the Theorem 3.1.1 (ignoring the genericity constraints {\bf 2}, {\bf 3} from Section 3.1) are invariant under $GL(n+2)$ acting on $\A^{n+2}$.

  {\bf Corollary 3.2.1.} $GL(n+2)$ acts on germs of admissible pairs at generic point.

  {\bf Definition 3.2.1.} A cone $C$ (as well as its projectivization) is called admissible if it is locally defined by (\ref{locs}) where function $f$ is admissible. The rank of $C$ is the rank of any such $f$. We denote the rank of $C$ by $rk(C)$. A cone $C$ is admissible iff $rk(C)>0$. We will also use notation $rk(\Sigma)$ for $rk(C)$ where $\Sigma\subset\P^{n+1}$ is the projectivization of $C$.

  The admissibility of cones defined above does not depend on the choice of projective coordinates.

  In order to write explicit formulas for $GL(n+2)$ action (see also (\ref{act})), it is convenient to describe $f,g$ parametrically as 
  \begin{equation}\label{par}
  f=\phi_0(u_1,...,u_n),~~~x_i=\phi_i(u_1,...,u_n),~i=1,...,n,
  \end{equation}
  $$g=\psi(u_1,...,u_n).$$
  Then the pair $\tilde{f}(x_1,...,x_n),\tilde{g}(x_1,...,x_n)$ defined parametrically by 
  \begin{equation}\label{par1}
  \tilde{f}=\tilde{\phi}_0(u_1,...,u_n),~~~x_i=\tilde{\phi}_i(u_1,...,u_n),~i=1,...,n,
  \end{equation}
  $$\tilde{g}=\tilde{\psi}(u_1,...,u_n)$$
  is also admissible, where
  \begin{equation}\label{actpr}
  \tilde{\phi}_i=\frac{\sum_{j=0}^na_{i,j}\phi_j+a_{i,n+1}}{\sum_{j=0}^na_{n+1,j}\phi_j+a_{n+1,n+1}},~i=0,...,n,
  \end{equation}
  $$\tilde{\psi}=\psi\cdot \frac{\det\Big(\frac{\partial \phi_i}{\partial u_j}\Big)_{1\leq i,j\leq n}}{\det\Big(\frac{\partial \tilde{\phi}_i}{\partial u_j}\Big)_{1\leq i,j\leq n}}\cdot\Big(\sum_{j=0}^na_{n+1,j}\phi_j+a_{n+1,n+1}\Big)^{-\frac{n+2}{2}}.$$
  Here $(a_{i,j})_{0\leq i,j\leq n+1}$ is an arbitrary non-degenerate constant matrix.

\

  {\bf Remark 3.2.1.} The projective invariance of admissible pairs is proven by a direct but not very transparent calculation. It can be explained in another way. First, notice that the group $Aff(n+1)$ of affine transformations of $\A^{n+1}$ acts on admissible pairs as
  \begin{equation}\label{actaf}
  \tilde{\phi}_i=\sum_{j=0}^na_{i,j}\phi_j+b_i,~i=0,...,n,
  \end{equation}
  $$\tilde{\psi}=\psi\cdot \frac{\det\Big(\frac{\partial \phi_i}{\partial u_j}\Big)_{1\leq i,j\leq n}}{\det\Big(\frac{\partial \tilde{\phi}_i}{\partial u_j}\Big)_{1\leq i,j\leq n}}.$$
  Here $(a_{i,j})_{0\leq i,j\leq n}$ is an arbitrary non-degenerate constant matrix and $b_i$ are arbitrary constants. 
  Notice that the last equation can be also written as\footnote{In invariant terms, we have a hypersurface in affine space $\A^{n+1}$, endowed with a volume element.} 
  $$\tilde{\psi}~ d\tilde{\phi}_1\wedge...d\tilde{\phi}_n=\psi~ d\phi_1\wedge...\wedge d\phi_n.$$
  Indeed, 
  after the change of variables 
  $$\hbar\mapsto \frac{\hbar}{y_0},~y_i\mapsto \frac{y_i}{y_0}$$
  the equation (\ref{eq0}) can be written as 
  \begin{equation}\label{eq2}
\int  e^{\frac{\phi_0 y_0+\phi_1 y_1+...+\phi_n y_n}{\hbar}}~\psi ~ d\phi_1\wedge...d\phi_n=(2\pi\hbar)^{\frac{n}{2}}\cdot y_0^{-\frac{n}{2}}~\hat{g}\Bigg(\frac{y_1}{y_0},...,\frac{y_n}{y_0}\Bigg)\cdot e^{\frac{y_0}{\hbar}\hat{f}\big(\frac{y_1}{y_0},...,\frac{y_n}{y_0}\big)}. 
  \end{equation}
  The l.h.s. of the equation (\ref{eq2}) is manifestly  invariant (up to multiplication by a constant independent of $u_1,...,u_n$) with respect to the affine action (\ref{actaf}) and the dual action of $GL(n+1)$ on variables $y_0,...,y_n$, so that the form $\phi_0 y_0+\phi_1 y_1+...+\phi_n y_n$ is invariant. 
  This gives the action (\ref{actaf}) of the group $Aff(n+1)$ on the set of admissible pairs. 

  By definition, the Fourier transform acts as an involution  (up to the reflection $x_i\to -x_i,i=1,\dots,n$) on the set of admissible pairs. Conjugating by the Fourier transform the action of the group $Aff(n+1)$ described above, we obtain the {\it second action} of the  same group  on the set of admissible pairs. One can check that these two actions generate the action of $GL(n+2)$.

  \subsection{Generalization to the projectively  dual lower-dimensional cones}

  It seems to be natural to generalize previous considerations in the following way:

  1. Densities $G(x_0,...,x_{n+1}),~\hat{G}(y_0,...,y_{n+1})$ are homogeneous, related by the Fourier transform but supported on projectively dual cones $C,\hat{C}$ of lower dimensions, not necessarily hypersurfaces.

  2. Densities $G,\hat{G}$ are both finite linear combinations of derivatives of delta functions, not necessarily just proportional to delta functions.

  We want to write explicitly conditions for $G,\hat{G}$. To simplify formulas, we assume that both $G,\hat{G}$ are proportional to delta functions, in the general case computations are similar.

  Suppose that our cone $C\subset\A^{n+2}$ is defined parametrically by
  $$x_i=\phi_i(u_0,...,u_{m_1}),~~~i=0,...,n+1$$
  where $u_0,...,u_{m_1}$ are coordinates on $C$, functions $\phi_i$ are homogeneous of degree $d_1$,  and $\dim C=m_1+1$.

  Similarly, suppose that  $\hat{C}\subset(\A^{n+2})^*$ is defined parametrically by
  $$y_i=\psi_i(v_0,...,v_{m_2}),~~~i=0,...,n+1$$
  where $v_0,...,v_{m_1}$ are coordinates on $\hat{C}$, functions $\psi_i$ are homogeneous of degree $d_2$,  and $\dim \hat{C}=m_2+1$.

  We can write our densities $G,\hat{G}$ as
  $$G(x_0,...,x_{n+1})=\int \prod_{i=0}^{n+1}\delta(x_i-\phi_i(u_0,...,u_{m_1}))g(u_0,...,u_{m_1})du_0...du_{m_1},$$
  $$\hat{G}(y_0,...,y_{n+1})=\int \prod_{i=0}^{n+1}\delta(y_i-\psi_i(v_0,...,v_{m_2}))\hat{g}(v_0,...,v_{m_2})dv_0...dv_{m_2}$$
  where $g,\hat{g}$ are also homogeneous with certain homogeneous degrees.
  The condition (\ref{FG}) after the integration with respect to variables  $x_0,...,x_{n+1}$ in the l.h.s. reads
  \begin{equation}\label{cond}
  \int g(u_0,...,u_{m_1})e^{\frac{1}{\hbar}\sum_{i=0}^{n+1}\phi_i(u_0,...,u_{m_1})y_i}du_0...du_{m_1}=
  \end{equation}
  $$(2\pi\hbar)^{\frac{n+2}{2}}\int \prod_{i=0}^{n+1}\delta(y_i-\psi_i(v_0,...,v_{m_2}))\hat{g}(v_0,...,v_{m_2})dv_0...dv_{m_2}.$$
  The homogeneity properties of the Fourier transform imply 
  $$\deg g+m_1+1=\frac{n+2}{2}~d_1,~~~\deg\hat{g}+m_2+1=\frac{n+2}{2}~d_2.$$
  One can deal with condition (\ref{cond}) in the following way. 

  a) Choose a partition $\{0,...,n+1\}=I_1\sqcup I_2$ such that $|I_1|=m_2+1$, such that $(y_i)_{i\in I_1}$ form a local system of coordinates on $\hat{C}$. 

  b) Perform integration in the r.h.s. with respect to $v_0,...,v_{m_2}$ removing delta functions $\delta(y_i-\psi_i(v_0,...,v_{m_2})),~i\in I_1.$ This reduces the number of delta functions in the r.h.s. by $m_2+1$.

  c) Multiply the equation (\ref{cond}) by $\frac{1}{(2\pi\hbar)^{n-m_2+1}}e^{\sum_{i\in I_2}\frac{y_iz_i}{\hbar}}$ and integrate with respect to $y_i,~i\in I_2$. After doing that the delta functions of the form $\prod_{i\in I_2}\delta(z_i+\phi_i(u_0,...,u_{m_1}))$ appear in the l.h.s., and we can remove them by integrating with respect to variables $u_j,~j\in J$ where  $J\subset \{0,...,m_1\}$ and $|J|=n-m_2+1$. In this way we remove all delta functions. Notice that we should have $n\leq m_1+m_2$, this is always true for projectively dual cones.

  d) After the removal of  all delta functions from (\ref{cond}) one can take an expansion of the l.h.s. at a critical point.

  {\bf Example 3.3.1.} Let $C\subset\A^5$ be defined parametrically by 
  $$C=\{(x_0,x_1,x_2,x_3.x_4)~\vert~x_2=h_1(u)x_0+h_2(u)x_1,x_3=h_3(u)x_0+h_4(u)x_1,x_4=h_5(u)x_0+h_6(u)x_1\}$$
  where $u$ is a coordinate on $C$ and functions $h_1,...,h_6$ satisfy conditions 
  $$h_4^{\prime}(u)h_1^{\prime}(u)=h_2^{\prime}(u)h_3^{\prime}(u),~~~h_6^{\prime}(u)h_1^{\prime}(u)=h_2^{\prime}(u)h_5^{\prime}(u)$$
  and certain genericity constraints. One can check that the dual cone $\hat{C}$ can also be defined parametrically as 
  $$\hat{C}=\{(y_0,y_1,y_2,y_3,y_4)~\vert~y_0=(p_1(v)h_1(v)-h_3(v))y_3+(p_2(v)h_1(v)-h_5(v))y_4,$$
  $$y_1=(p_1(v)h_2(v)-h_4(v))y_3+(p_2(v)h_2(v)-h_6(v))y_4,y_2=-p_1(v)y_3-p_2(v)y_4\}$$
  where $v$ is a coordinate on $\hat{C}$ and 
  $$p_1(v)=\frac{h_3^{\prime}(v)}{h_1^{\prime}(v)}=\frac{h_4^{\prime}(v)}{h_2^{\prime}(v)},~~~p_2(v)=\frac{h_5^{\prime}(v)}{h_1^{\prime}(v)}=\frac{h_6^{\prime}(v)}{h_2^{\prime}(v)}.$$
  We have $\dim C=\dim \hat{C}=3$ and therefore both $C$ and its dual $\hat{C}$ have codimension 2 in $\A^5$.

  Admissibility condition (\ref{cond}) in this case can be written as 
  $$\int g(x_0,x_1,u)e^{\frac{1}{\hbar}(x_0y_0+x_1y_1+(h_1(u)x_0+h_2(u)x_1)y_2+(h_3(u)x_0+h_4(u)x_1)y_3+(h_5(u)x_0+h_6(u)x_1)y_4}dx_0dx_1du$$
  $$=(2\pi\hbar)^{\frac{5}{2}}\int\hat{g}(y_3,y_4,v)\delta(y_2+p_1(v)y_3+p_2(v)y_4)$$
  $$\delta((p_1(v)h_1(v)-h_3(v))y_3+(p_2(v)h_1(v)-h_5(v))y_4-y_0)$$
  $$\delta((p_1(v)h_2(v)-h_4(v))y_3+(p_2(v)h_2(v)-h_6(v))y_4-y_1)dv.$$
  Multiplying both sides by $\frac{1}{(2\pi\hbar)^2}e^{-\frac{1}{\hbar}(y_0z_0+y_1z_1)}$, integrating by $y_0,y_1$, and in the l.h.s. integrating also by $x_0,x_1$ we obtain
  $$\int g(z_0,z_1,u)e^{\frac{1}{\hbar}((h_1(u)z_0+h_2(u)z_1)y_2+(h_3(u)z_0+h_4(u)z_1)y_3+(h_5(u)z_0+h_6(u)z_1)y_4)}du=$$
  $$(2\pi\hbar)^{\frac{1}{2}}\int\hat{g}(y_3,y_4,v)\delta(y_2+p_1(v)y_3+p_2(v)y_4)$$
  $$e^{-\frac{1}{\hbar}(((p_1(v)h_1(v)-h_3(v))y_3+(p_2(v)h_1(v)-h_5(v))y_4)z_0+((p_1(v)h_2(v)-h_4(v))y_3+(p_2(v)h_2(v)-h_6(v))y_4)z_1)}dv.$$
  We can integrate by $v$ in the r.h.s. removing delta function, and integrate by $u$ in the l.h.s. by taking expansion at a critical point. In this way we obtain admissibility conditions for $C$ explicitly. It would be interesting to classify such admissible cones and generalize these for higher dimensions.

  \section{Reformulation of the problem in terms of constraints on wave functions}

  \subsection{Abstract formalism}

  In the previous Section we reformulated our original problem of finding admissible pairs $f,g$ in $n$ variables in terms of finding distributions $G,\hat{G}$ in $n+2$ variables, given by (\ref{G}), (\ref{Ghat}), which are related by the Fourier transform (\ref{FG}). Let us apply the formalism of wave functions from Section 2 to these distributions. 

  Conical germs $C$ and $\hat{C}$ on which distributions  $G$, $\hat{G}$ are supported, are given by equations 
  \begin{equation}\label{eqC}
  \frac{x_{n+1}}{x_0}-f\Big(\frac{x_1}{x_0},...,\frac{x_n}{x_0}\Big)=0,~~~~~~~\frac{y_0}{y_{n+1}}+\hat{f}\Big(\frac{y_1}{y_{n+1}},...,\frac{y_n}{y_{n+1}}\Big)=0.
  \end{equation}

  Denote by $L$ the conormal bundle to cone $C\subset\A^{n+2}$. It is a Lagrangian subvariety of $\A^{2(n+2)}=T^*\A^{n+2}$ invariant under the action of $GL(1)\times GL(1)$ on $\A^{2(n+2)}$ where the first $GL(1)$ rescales coordinates $x_0,...,x_{n+1}$ and the second $GL(1)$ rescales coordinates $y_0,...,y_{n+1}$. In particular, $L$ is conical\footnote{Generically, the invariance under the second copy of $GL(1)$ means that $L$ is a conormal bundle to a subvariety in $\A^{n+2}$ of arbitrary dimension. The invariance under the first copy of $GL(1)$ means that this subvariety is conical.}.

  Nonvanishing of the Hessian of germ $f$ (and hence of $\hat{f}$) can be reformulated as the following genericity condition on the smooth $GL(1)\times GL(1)$ invariant germ $L$ of Lagrangian submanifold at point $p\in\A^{2(n+2)}$, similar to conditions {\bf 1, 2, 3} in Section 3.1:

  1) coordinates $x_0$ and $y_{n+1}$ of $p$ are non-zero. 

  2) projections from $T_pL$ to $\A^{n+2}_{x_0,...,x_n,y_{n+1}}$ and to $\A^{n+2}_{y_0,...,y_n,x_{n+1}}$ are one-to-one.

  The space $\A^{2(n+2)}$ with the symplectic structure $\sum_{i=0}^{n+1} dx_i\wedge dy_i$ is simultaneously the cotangent space to $\A^{n+2}_{x_0,...,x_{n+1}}$ and to $\A^{n+2}_{y_0,...,y_{n+1}}$. 

  Projective duality between the projectivizations of cones $C$ and $\hat{C}$ can be reformulated as the property of $L$ to be simultaneously the conormal bundle to $C$ and $\hat{C}$. 

  Let  
  \begin{equation}\label{eqC1}
  F_1(x_0,...,x_{n+1})=0, ~~~~~~ F_2(y_0,...,y_{n+1})=0 
  \end{equation}
  be equations 
  for the cones $C$ and $\hat{C}$ respectively. Here $F_1,~F_2$ are germs at non-zero points of homogeneous functions of some homogeneous degrees which are proportional with an invertible factors to l.h.s. of equations (\ref{eqC}).

  Distributions $G,\hat{G}$ can be understood as elements $\psi, \hat{\psi}\in {\mathcal WF}_{L}$. Similarly to Remark 2.2 the fact that $\hat{G}$ is the Fourier transform of $G$ can be written as equality $\psi=\hat{\psi}$. 

  Let us interpret $F_1,~F_2$ as elements of the quantum algebra acting on ${\mathcal WF}_{L}$. The property of $G$ to be a smooth density on $C$ can be rewritten as\footnote{Recall that $C$ is defined by equation $F_1=0$ and therefore, $G$ is proportional to $\delta(F_1)$. For example, we have $P(x_0,...,x_{n+1})\delta(P(x_0,...,x_{n+1}))=0$ for an arbitrary polynomial $P$.} 
  
  {\bf 1.} $F_1\cdot \psi=0.$

  Similarly the dual condition on $\hat{G}$ gives
  
  ${\bf 1^{\prime}.}$  $F_2\cdot\psi=0.$

  The condition on $G,\hat{G}$ of being homogeneous of degree $-\frac{n+2}{2}$ can be rewritten as

  {\bf 2.} $\Big(\sum_{i=0}^{n+1}y_i\star x_i\Big)\cdot \psi=0.$

  Indeed, the homogeneity of $\psi$ (considered as a distribution in variables $x_0,...,x_{n+1}$) of degree $-\frac{n+2}{2}$ means 
  $$\sum_{i=0}^{n+1}x_i\partial_{x_i}\psi+\frac{n+2}{2}\psi=0,$$ 
  therefore 
  $\Big(\sum_{i=0}^{n+1}y_i\star x_i\Big)\cdot \psi=\frac{1}{2}\Big(\sum_{i=0}^{n+1}\hbar\partial_{x_i} x_i+\hbar x_i\partial_{x_i}\Big)\cdot \psi=\frac{\hbar}{2}\sum_{i=0}^{n+1}\Big(2x_i\partial_{x_i}+1\Big)\psi=0.$

  Finally, the property that $G,\hat{G}$ are given by densities independent of $\hbar$ can be rewritten as 

  {\bf 3.} $\tau_L\psi=-\frac{n}{4}\psi$ where $\tau_L$ is defined in Remark 2.3 in Section 2. The eigenvalue $-\frac{n}{4}$ can be seen from the formula (\ref{fe}) in Section 4.2.

  \ 

  The discussion above can be summarized as 

  {\bf Proposition 4.1.1.} Fix functions $f,\hat{f}$ and the corresponding cones $C,\hat{C}$ given by equations (\ref{eqC1}). There is one-to-one correspondence between the space of $g$ such that pair $f,g$ is admissible and non-zero elements $\psi \in {\mathcal WF}_{L}$ satisfying properties ${\bf 1,~1^{\prime}~,2,~3}$ where $L$ is a $GL(1)\times GL(1)$-invariant Lagrangian germ satisfying the genericity condition above. $\square$

  \subsection{Explicit formulas in general case}

  Let us write equations on $\psi$ explicitly.  Conditions {\bf 1, 2, 3} give
  \begin{equation}\label{psi}
  \psi=g\Big(\frac{x_1}{x_0},...,\frac{x_n}{x_0}\Big)~x_0^{-\frac{n+2}{2}}\delta\Big(\frac{x_{n+1}}{x_0}-f\Big(\frac{x_1}{x_0},...,\frac{x_n}{x_0}\Big)\Big)
  \end{equation}
  for some function $g$. Indeed, condition {\bf 1} means that $\psi$ is proportional to the delta function of cone $C$, condition {\bf 2} gives homogeneity condition for the coefficient of proportionality, and condition {\bf 3} means that function $g$ does not depend on $\hbar$. However, condition ${\bf 1^{\prime}}$ becomes singular in these coordinates because the action of $F_2(\hbar\partial_{x_0},...,\hbar\partial_{x_{n+1}})$ in general is not defined on $\psi$ given by (\ref{psi}). In order to overcome this problem, let us make the Fourier transform with respect to the coordinate $x_{n+1}$. Geometrically this means that we choose the projection of our cone $C$ to coordinates $x_0,...,x_n,y_{n+1}$.
  After making the Fourier transform in $x_{n+1}$ the element $\psi$ given by (\ref{psi}) takes the form 
   \begin{equation}\label{fe}
  \psi=x_0^{-\frac{n}{2}}e^{\frac{x_0y_{n+1}}{\hbar}f\big(\frac{x_1}{x_0},...,\frac{x_n}{x_0}\big)}g\Big(\frac{x_1}{x_0},...,\frac{x_n}{x_0}\Big)
  \end{equation}
  and condition ${\bf 1^{\prime}}$ means that
  $$F_2(\hbar \partial_{x_0},...\hbar\partial_{x_n},y_{n+1}) \cdot \psi=0$$
  and can be rewritten as an equation on $g$
  \begin{equation}\label{eqgen}
  F_2\Big(\hbar\partial_{x_0}-\frac{n\hbar}{2x_0}+Q_0,\hbar\partial_{x_1}+Q_1,...,\hbar\partial_{x_n}+Q_n,y_{n+1}\Big)\cdot g\Big(\frac{x_1}{x_0},...,\frac{x_n}{x_0}\Big)=0
  \end{equation}
  where operator $F_2$ in (\ref{eqgen}) is equal to
  $$\Bigg(x_0^{-\frac{n}{2}}e^{\frac{x_0y_{n+1}}{\hbar}f\big(\frac{x_1}{x_0},...,\frac{x_n}{x_0}\big)}\Bigg)^{-1} \cdot F_2(\hbar \partial_{x_0},...\hbar\partial_{x_n},y_{n+1})\cdot x_0^{-\frac{n}{2}}e^{\frac{x_0y_{n+1}}{\hbar}f\big(\frac{x_1}{x_0},...,\frac{x_n}{x_0}\big)}$$
  and, therefore, $Q_0,...,Q_n$ are  defined by 
  $$\Bigg(x_0^{-\frac{n}{2}}e^{\frac{x_0y_{n+1}}{\hbar}f\big(\frac{x_1}{x_0},...,\frac{x_n}{x_0}\big)}\Bigg)^{-1}\cdot \hbar \partial_{x_0}\cdot x_0^{-\frac{n}{2}}e^{\frac{x_0y_{n+1}}{\hbar}f\big(\frac{x_1}{x_0},...,\frac{x_n}{x_0}\big)}=\hbar\partial_{x_0}-\frac{n\hbar}{2x_0}+Q_0,$$
  $$\Bigg(x_0^{-\frac{n}{2}}e^{\frac{x_0y_{n+1}}{\hbar}f\big(\frac{x_1}{x_0},...,\frac{x_n}{x_0}\big)}\Bigg)^{-1}\cdot \hbar \partial_{x_i}\cdot x_0^{-\frac{n}{2}}e^{\frac{x_0y_{n+1}}{\hbar}f\big(\frac{x_1}{x_0},...,\frac{x_n}{x_0}\big)}=\hbar\partial_{x_i}+Q_i,~~~i=1,...,n.$$
  Explicitly, we have
  $$Q_0=y_{n+1}f\Big(\frac{x_1}{x_0},...,\frac{x_n}{x_0}\Big)-\frac{y_{n+1}}{x_0}\sum_{i=1}^n x_if_i\Big(\frac{x_1}{x_0},...,\frac{x_n}{x_0}\Big),$$
  $$Q_i=y_{n+1}f_i\Big(\frac{x_1}{x_0},...,\frac{x_n}{x_0}\Big),~~~i=1,...,n$$
  where by $f_i$ we denote partial derivative of $f$ with respect to its $i-$th argument.
  
  Notice that $Q_0,...,Q_n$ are independent of $\hbar$. 

  The l.h.s. of (\ref{eqgen}) makes sense as a power series in $\hbar$. For example, if we choose $F_2$ in the simplest form 
  $$F_2(y_0,...,y_{n+1})=\frac{y_0}{y_{n+1}}+\hat{f}\Big(\frac{y_1}{y_{n+1}},...,\frac{y_n}{y_{n+1}}\Big)$$
  as in (\ref{eqC}), and take Taylor expansion\footnote{Computation of this Taylor series is not a straightforward problem because while arguments of $F_2$ pairwise commute, parts of these arguments at different powers of $\hbar$ does not, for example $\hbar\partial_{x_i}$ does not commute with $Q_i$.} of the l.h.s. of (\ref{eqgen}) at $\hbar=0$ we obtain a power series in $\hbar$ where each term is a differential operator (with coefficients written in terms of $f$) applied to the function $g$. This can be written as a system of partial differential equations for unknown functions $f,g$ linear in $g$. Notice that this system is infinite in general. 

  \subsection{Algebraic case}

  Consider the {\it algebraic} case where $f$ is an algebraic function and, therefore, $\hat{f},~C,~\hat{C}$ are algebraic. Let us choose $F_1(x_0,...,x_{n+1}),~F_2(y_0,...,y_{n+1})$ to be homogeneous polynomials of certain degrees, such that 
  $$F_1(1,x_1,...,x_n,f(x_1,...,x_n))=0,~~~~~~~F_2(1,y_1,...,y_n,\hat{f}(y_1,...,y_n))=0.$$

  In this special case one can eliminate parameter $\hbar$, and the problem of finding an element $\psi$ satisfying to properties ${\bf 1,~1^{\prime},~2,~3}$ turns to the following purely algebraic question. 

  Find homogeneous polynomials $F_1(x_0,...,x_{n+1}),~F_2(y_0,...,y_{n+1})$ in two dual groups of variables (or more abstractly, in dual vector spaces) such that there exists a non-zero cyclic module $M$ over the ring of polynomial differential operators 
  $${\bf k}[x_0,...,x_{n+1}][\partial_{x_0},...,\partial_{x_{n+1}}]$$
  generated by an element $\psi\in M$ such that 
  \begin{equation}
  \label{eqalg}
  \begin{split}
  &F_1(x_0,...,x_{n+1})\cdot(\psi)=0,\\
   &F_2(\partial_{x_0},...,\partial_{x_{n+1}})\cdot(\psi)=0,\\
  &\Big(\sum_{i=0}^{n+1}x_i\partial_{x_i}+\frac{n+2}{2}\Big)\psi=0.
  \end{split}
  \end{equation}
  Notice that the first equation means that 
  $$\psi=G_1\Big(\frac{x_1}{x_0},...,\frac{x_n}{x_0}\Big)x_0^m\delta(F_1(x_0,...,x_{n+1}))$$
  where the factor in front of the delta function is a homogeneous function of degree $m$. The third equation means that 
  $$m-\deg F_1=-\frac{n+2}{2}\,.$$
  The second equation gives a finite system of partial differential equations on $F_1,G_1$ linear in $G_1$.

  One has a universal module $M_{\psi}$ generated by the cyclic vector $\psi$ satisfying the above equations.

  Roughly speaking, the study of admissible pairs with algebraic function $f$ (or equivalently, $\hat{f}$) can be rephrased as the study of projectively dual algebraic hypersurfaces in $\P^{n+1}$, $(\P^{n+1})^*$ such that the corresponding finitely generated module $M_{\psi}$ is not zero. 

  We will see in the sequel both holonomic and nonholonomic\footnote{Holonomicity at generic point means that the space of $G_1$ is finite-dimensional.} examples of modules $M_{\psi}$.

\

  {\bf Remark 4.3.1.} Let $\Sigma\subset\P^{n+1}$ and $\widehat{\Sigma}\subset(\P^{n+1})^*$ be projectively dual hypersurfaces. Let $\Sigma$ (resp, $\widehat{\Sigma}$) be defined by a homogeneous irreducible polynomial $F_1(x_0,...,x_{n+1})$ (resp. $F_2(y_0,...,y_{n+1})$). Given polynomial $F_1$ one can determine polynomial $F_2$ (up to  multiplication by a non-zero constant) from the condition: $F_2$ is a homogeneous polynomial of smallest degree such that 
  $$F_2(\partial_{x_0}F_1,...,\partial_{x_{n+1}}F_1)=0 \mod F_1,$$
  and similarly one can find polynomial $F_1$ if polynomial $F_2$ is given. We can define polynomials $H_1\in {\bf k}[x_0,...,x_{n+1}]$ and $H_2\in {\bf k}[y_0,...,y_{n+1}]$ by 
   \begin{equation}\label{hh}
  F_1(\partial_{y_0}F_2,...,\partial_{y_{n+1}}F_2)=(f_2-1)F_2(y_0,...,y_{n+1})H_2(y_0,...,y_{n+1}),
  \end{equation}
  $$F_2(\partial_{x_0}F_1,...,\partial_{x_{n+1}}F_1)=(f_1-1)F_1(x_0,...,x_{n+1})H_1(x_0,...,x_{n+1})$$
  where $f_1=\deg F_1,~f_2=\deg F_2$.
   Notice that 
  $$\deg H_1=\deg H_2=f_1f_2-f_1-f_2.$$

   The mappings 
  $$x_i\mapsto \partial_{y_i}F_2,~~~y_i\mapsto \partial_{x_i} F_1,~~~i=1,...,n$$
  define birational isomorphisms between projective hypersurfaces $\Sigma,~\widehat{\Sigma}$. They are mutually inverse as birational mappings of projective varieties. This can be written algebraically as 
  $$\partial_{y_i}F_2~\Big\rvert_{y_0=\partial_{x_0}F_1,...,y_{n+1}=\partial_{x_{n+1}}F_1}=x_iH_1~\mod~ F_1,~~~i=1,...,n$$
  and
  $$\partial_{x_i}F_1~\Big\rvert_{x_0=\partial_{y_0}F_2,...,x_{n+1}=\partial_{y_{n+1}}F_2}=y_iH_2~\mod~ F_2,~~~i=1,...,n.$$

  These formulas can be verified by computing derivatives of equations (\ref{hh}), reducing modulo $F_1$ or $F_2$ and using Euler's homogeneous function theorem.

  \ 

  {\bf Remark 4.3.2.} Projective duality interchanges polynomials $F_1(x_0,...,x_{n+1}),~F_2(y_0,...,y_{n+1})$. Dual element $\hat{\psi}$ has a form (up to a constant) 
  $$\hat{\psi}=G_1\Big(\frac{x_1}{x_0},...,\frac{x_n}{x_0}\Big)x_0^m~\Big\rvert_{x_0=\partial_{y_0}F_2,...,x_{n+1}=\partial_{y_{n+1}}F_2}\det\Big(\partial_{y_i}\partial_{y_j}F_2\Big)^{\frac{1}{2}}H_2^{-1}\delta(F_2(y_0,...,y_{n+1}))$$
  where $H_2\in {\bf k}[y_1,...,y_n]$ (and similar polynomial $H_1$) are defined above in Remark 4.3.1.

  \ 

  {\bf Remark 4.3.3.} The question of finding pairs of polynomials 
  $$F_1\in {\bf k}[x_0,...,x_{n+1}],~~~F_2\in {\bf k}[\partial_{x_0},...,\partial_{x_{n+1}}]$$ generating a non-zero cyclic module $M$ over the ring of differential operators, can be generalized to an associative algebra with two maximal commutative subalgebras. In the case of differential operators homogeneity of $F_1,~F_2$ follows automatically from the nonvanishing of $M$ (assuming that both $F_1,~F_2$ are irredusible). In other natural cases, such as quantum torus or ring of difference operators we are not aware of any nontrivial example.

\ 

  {\bf Remark 4.3.4.} The condition that both $F_1,~F_2$ are irreducible polynomials defining projectively dual non-degenerate hypersurfaces can be generalized. 

  First, one can consider the case where $F_1$ is a power of irreducible polynomial. In terms of distributions $G,\hat{G}$ it means that they are finite linear combinations of derivatives of delta functions supported on corresponding coves. In the original formulation it means that $g,\hat{g}$ are polynomials in $\hbar$.

  Second, for a pair of projectively dual projective manifolds one (or even both) of them could have codimension larger than one, see Section 3.3.

  \subsection{Explicit formulas in the algebraic case} 

  Let us choose $F_1$ in the non-polynomial form
  $$F_1=x_{n+1}-x_0f\Big(\frac{x_1}{x_0},...,\frac{x_n}{x_0}\Big)$$
  where $f$ is an algebraic function, and 
  $$F_2\in {\bf k}[\partial_{x_0},...,\partial_{x_{n+1}}]$$
  is an irreducible polynomials as above. In this case properties {\bf 1,~2,~3} give (see (\ref{psi}))
  \begin{equation}\label{psi1}
  \psi=g\Big(\frac{x_1}{x_0},...,\frac{x_n}{x_0}\Big)~x_0^{-\frac{n}{2}}\delta\Big(x_{n+1}-x_0f\Big(\frac{x_1}{x_0},...,\frac{x_n}{x_0}\Big)\Big)
  \end{equation}
  and the property ${\bf 1^{\prime}}$ reads
  \begin{equation}\label{eqpol}
  F_2(\partial_{x_0},...,\partial_{x_{n+1}})\cdot(\psi)=0
  \end{equation}
  where $\psi$ is given by (\ref{psi1}). Cyclic $D$-module over the ring ${\bf k}[\partial_{x_0},...,\partial_{x_{n+1}}]$ generated by $\psi$ is contained in the bigraded vector space whose $(a,b)-$graded component is spent by expressions of the form 
  $$h\Big(\frac{x_1}{x_0},...,\frac{x_n}{x_0}\Big)~x_0^{a}\delta^{(b)}\Big(x_{n+1}-x_0f\Big(\frac{x_1}{x_0},...,\frac{x_n}{x_0}\Big)\Big)$$
  where $h$ is a function in $n$ variables, $a\in -\frac{n}{2}-\Z_{\geq 0}$,  $b\in\Z_{\geq0}$
  We can identify the above bigraded vector space as the module over ${\bf k}[\partial_{x_0},...,\partial_{x_{n+1}}]$ with another one spent by 
  $$h\Big(\frac{x_1}{x_0},...,\frac{x_n}{x_0}\Big)~\frac{x_0^{a}}{\Big(x_{n+1}-x_0f\Big(\frac{x_1}{x_0},...,\frac{x_n}{x_0}\Big)\Big)^{b+1}}$$
  where $h,a,b$ are as above. The isomorphism is given by $s\delta^{(b)}(t)\mapsto s\frac{(-1)^bb!}{t^{b+1}}$.

  The l.h.s. of (\ref{eqpol}) lies in the direct sum of bigraded components with $a=-\frac{n}{2}-k,~b=\deg F_2-k$ where $k=0,...,\deg F_2$.

  Introduce new variables $v_0,...,v_{n+1}$ such that
  $$v_i=\frac{x_i}{x_0},~i=1,...,n,~~~~v_0=x_0,~~~~v_{n+1}=x_{n+1}-x_0f\Big(\frac{x_1}{x_0},...,\frac{x_n}{x_0}\Big).$$
  In these variables element $\psi$ given by (\ref{psi1}) reads as
  \begin{equation}\label{psi2}
  \psi=g\big(v_1,...,v_n\big)~v_0^{-\frac{n}{2}}v_{n+1}^{-1}.
  \end{equation}
  Partial derivatives $\partial_{x_0},...,\partial_{x_{n+1}}$ in new variables become more complicated differential operators 
  $$\partial_{x_i}\mapsto D_i,~i=0,...,n+1$$ 
  where 
  $$D_0=-\frac{1}{v_0}\sum_{j=1}^nv_j\partial_{v_j}+\partial_{v_0}+\Big(\sum_{j=1}^nv_jf_{v_j}-f\Big)\partial_{v_{n+1}}$$
  $$D_i=\frac{1}{v_0}\partial_{v_i}-f_{v_i}\partial_{v_{n+1}},~i=1,...,n,~~~~D_{n+1}=\partial_{v_{n+1}}$$
  where $f=f(v_1,...,v_n),~f_{v_i}=\frac{\partial f}{\partial v_i}$.

  The final conclusion is that in coordinates $v_0,...,v_{n+1}$ the system of equations for $f,g$ can be written as 
  \begin{equation}\label{eq}
  F_2(D_0,...,D_{n+1})\cdot \Big(g(v_1,...,v_n)v_0^{
  -\frac{n}{2}}v_{n+1}^{-1}\Big)=0.
  \end{equation}
  To obtain a system of partial differential equations for $f,g$ one should equate coefficients at all powers of $v_0,v_{n+1}$ in (\ref{eq}) to zero. In this way we get a system of $\deg F_2$ differential equations on $f,g$. 

  Notice that if our cone $C$ is given parametrically by (\ref{par}), we can rewrite our system (\ref{eq}) in parametric form by doing the corresponding change of variables. 
  
  \subsection{A question about monodromic regular holonomic $D$-modules}

  One can try to look for admissible pairs of a special type, when both homogeneous distributions $G,~\hat{G}$ generate {\it regular holonomic} $D$-modules.

  Recall the classical result by J. L. Brylinski \cite{bk} which says that regular holonomic $D$-modules on a vector space $V$ such that their Fourier transforms are also regular holonomic are exactly those for which the action of Euler vector field $E=\sum_ix_i\partial_{x_i}$ is locally finite. Such $D$-modules are called {\it monodromic regular holonomic}. 

  The singular support $SS(M)$ of such a $D$-module $M$ is $GL(1)\times GL(1)$-invariant (possibly reducible) Lagrangian cone $L\subset T^*V=V\oplus V^*$, and it coincides with the singular support of its Fourier transform $SS({\cal F}(M))$ under the identification $T^*V=V\oplus V^*=T^*V^*$.

  Thus, we arrive to the following problem: study monodromic regular holonomic $D$-modules $M$ such that $SS(M)$ does not contain $V\times\{0\}$ and $\{0\}\times V^*$. Indeed, in this case any non-zero element $\psi\in M$ homogeneous with respect to $E$ is killed by some non-trivial homogeneous polynomial $F_1(x_i),~F_2(\partial_{x_i})$.

  Riemann-Hilbert correspondence identifies monodromic regular holonomic $D$-modules with so-called monodromic perverse sheaves on $V$. Thus, one can reformulate the above problem in purely topological terms, concerning finite-dimensional representations of the fundamental group of $L\setminus L^{sing}$. 

  Finally, for an irreducible algebraic cone $C\subset V$, such that $C\ne 0,V$ there is a natural $GL(1)\times GL(1)$ invariant irreducible Lagrangian cone $L$ which is the conormal bundle to $C$. Notice that $L=L_C\ne V\times\{0\},~\{0\}\times V^*$. So $L_C$ is a natural candidate for the singular support.

  {\bf Question 4.5.1.} For which $C$ there exists a monodromic regular holonomic $D$-module whose singular support is $L_C$, may be with multiplicities?

 \section{Examples of admissible hypersurfaces and corresponding pairs}

 In this Section we assume that $\bf k=\bar{k}$.
 
 \subsection{Quadratic hypersurfaces}
 
 {\bf Lemma 5.1.1.} Let $f(x_1,...,x_n)=-\frac{1}{2}(x_1^2+...+x_n^2)$. Then pair $f,g$ is admissible iff the function $g(x_1,...,x_n)$ is harmonic, i.e.
 $$\frac{\partial^2g}{\partial x_1^2}+...+\frac{\partial^2g}{\partial x_n^2}=0.$$
 We have in this case 
 $$\int g(x_1,...,x_n)e^{\frac{1}{\hbar}(-\frac{1}{2}x_1^2-...-\frac{1}{2}x_n^2+x_1y_1+...+x_ny_n)}dx_1...dx_n=(2\pi \hbar)^{\frac{n}{2}}g(y_1,...,y_n)e^{\frac{1}{2\hbar}(y_1^2+...+y_n^2)}.$$
 
 {\bf Proof.} For arbitrary function $g$ we have
 $$\int g(x_1,...,x_n)e^{\frac{1}{\hbar}(-\frac{1}{2}x_1^2-...-\frac{1}{2}x_n^2+x_1y_1+...+x_ny_n)}dx_1...dx_n=(2\pi \hbar)^{\frac{n}{2}}e^{\frac{1}{2\hbar}(y_1^2+...+y_n^2)}\sum_{i=0}^{\infty}\frac{\hbar^i}{i!}\Delta^ig(y_1,...,y_n)$$
 where $\Delta=\frac{\partial^2}{\partial y_1^2}+...+\frac{\partial^2}{\partial y_n^2}$. The r.h.s. consists only of the first term iff $\Delta g=0$.
 $\square$
 
 {\bf Theorem 5.1.1.}  Let $Q=\sum_{i,j=0}^{n+1}a_{i,j}x_ix_j$ be a non-degenerate quadric. Here $a_{i,j}=a_{j,i}$ and $\det(a_{i,j})\ne 0$. Then the hypersurface in $\P^{n+1}$ 
 defined by $Q=0$ is admissible. Its rank is two if $n=1$ and infinity otherwise.
 
 {\bf Proof.} Any such quadric is $GL(n+2)$ equivalent to $Q_0=2x_0x_{n+1}+x_1^2+...+x_n^2$. The projective hypersurface $Q_0=0$ is a projectivization of affine hypersurface
 defined by the equation $x_{n+1}=-\frac{1}{2}(x_1^2+...+x_n^2)$ which is admissible by Lemma 1. Therefore, the hypersurface $Q=0$ is also admissible by projective invariance.  $\square$

 \ 
 
 {\bf Remark 5.1.1.} Solving the equation 
 $$Q(1,x_1,...,x_n,f)=0$$
 for an arbitrary non-degenerate quadric $Q$ we obtain a family of functions $f$ which depends on a lot of parameters $a_{i,j}$. One can see that by applying group of transformations 
 $$x_i\mapsto \sum_{j=1}^nq_{i,j}x_j+b_j,~~~f\mapsto \lambda f+\sum_{j=1}^n\mu_jx_j+\nu$$
 one can reduce $f$ to one of the following forms:
 
 {\bf a)}  $f(x_1,...,x_n)=(x_1^2+...+x_n^2+1)^{\frac{1}{2}}$ ,
 
 {\bf b)}  $f(x_1,...,x_n)=(x_1^2+...+x_{n-1}^2+x_n)^{\frac{1}{2}}$  for $n>1$, and  $f(x_1)=x_1^{\frac{1}{2}}$ for $n=1$,
 
 {\bf c)}  $f(x_1,...,x_n)=\frac{x_1^2+...+x_{n-1}^2+1}{x_n}$ for $n>1$, and $f(x_1)=\frac{1}{x_1}$ for $n=1$,
 
 {\bf d)}  $f(x_1,...,x_n)=-x_1^2-...-x_n^2$.

 In each of the cases above pair $f,g$ is admissible iff $g$ satisfies a certain second order partial differential equation which can be obtained from harmonic equation in Lemma 5.1.1 by the action of $GL(n+2)$, see Section 3.2.

\ 

 It would be interesting to lift, if possible, the corresponding formal integrals of the form (\ref{eq0}) to actual convergent integral identities. For example, it is known \cite{b} that if
 $$ g(\vec{x})=\frac{1}{\sqrt{1+\vec x^2}\cdot (1+\sqrt{1+\vec{x}^2})^{\frac{n-2}{2}}}$$
 then
 $$\int_{\R^n}g(\vec{x})e^{-\frac{1}{\hbar}\sqrt{1+\vec{x}^2}+\frac{i}{\hbar}\vec{x}\vec{y}}d\vec{x}=(2\pi \hbar)^{\frac{n}{2}}g(\vec{y})e^{-\frac{1}{\hbar}\sqrt{1+\vec{y}^2}}$$
 
  \subsection{Hypersurfaces admitting quadratic parametrization}

  Let  $C\subset\A^{n+2}$ be a cone given parametrically by
  $$x_i=\frac{1}{2}\sum_{j,k=0}^na_{i,j,k}u_ju_k,~~~i=0,...,n+1$$
  where $a_{i,j,k}\in {\bf k}$ are constants and $a_{i,k,j}=a_{i,j,k}.$
  The dual cone $\hat{C}$ is given by 
  $$\det\Bigg(\sum_{i=0}^{n+1}a_{i,j,k}y_i\Bigg)_{0\leq j,k\leq n}=0.$$
  Indeed, let $Q=x_0y_0+...+x_{n+1}y_{n+1}=\frac{1}{2}\sum_{i=0}^{n+1}\sum_{j,k=0}^na_{i,j,k}y_iu_ju_k.$ The point $(y_0,...,y_{n+1})\in\A^{n+2}$ belongs to the dual cone $\hat{C}$ iff the linear system $\partial_{u_i}Q=0,~i=0,...,n$ for $u_0,...,u_n$ has a non-zero solution.

  It follows from Theorem 3.1.1 that cone $C$ is admissible iff
  $$\int g_1\Big(\frac{x_1}{x_0},...,\frac{x_n}{x_0}\Big)~x_0^{-\frac{n}{2}}e^{\frac{1}{\hbar}(x_0y_{n+1}f\big(\frac{x_1}{x_0},...,\frac{x_n}{x_0}\big)+x_0y_0+...+x_{n}y_{n})}dx_0...dx_n=$$
  $$(2\pi\hbar)^{\frac{n+2}{2}}\delta\Big(\frac{y_0}{y_{n+1}}+\hat{f}\Big(\frac{y_1}{y_{n+1}},...,\frac{y_n}{y_{n+1}}\Big)\Big)g_2\Big(\frac{y_1}{y_{n+1}},...,\frac{y_n}{y_{n+1}}\Big)y_{n+1}^{-\frac{n+2}{2}}$$
  for some functions $g_1,g_2$, where $x_0y_{n+1}f\big(\frac{x_1}{x_0},...,\frac{x_n}{x_0}\big)=x_{n+1}$ and $\hat{f}$ corresponds to the dual cone $\hat{C}$. Substituting our quadratic parametrization in the l.h.s. and taking into account the equation for the dual cone we get
  $$\int g(u_0,...,u_n)e^{\frac{1}{2\hbar}\sum_{i=0}^{n+1}\sum_{j,k=0}^na_{i,j,k}y_iu_ju_k}du_0...du_n=$$
  $$(2\pi\hbar)^{\frac{n+2}{2}}\delta\Bigg(\det\Bigg(\sum_{i=0}^{n+1}a_{i,j,k}y_i\Bigg)_{0\leq j,k\leq n}\Bigg)\hat{g}\Big(\frac{y_1}{y_{n+1}},...,\frac{y_n}{y_{n+1}}\Big)y_{n+1}^{\frac{n}{2}}$$
  where $g$ is a homogeneous function of degree 1.

  {\bf Theorem 5.2.1.} Any non-degenerate cone $C$ admitting a quadratic parametrization is admissible with the rank larger or equal than $n+1$. Moreover, any linear homogeneous function $g$ is admissible.

  {\bf Proof.} Without any loss of generality we can set $g(u_0,...,u_n)=u_0$ because of the $GL(n+1)$-action on variables $u_0,...,u_n$. In order to compute the integral 
  \begin{equation}\label{int}
   \int u_0e^{\frac{1}{2\hbar}\sum_{i=0}^{n+1}\sum_{j,k=0}^na_{i,j,k}y_iu_ju_k}du_0...du_n
  \end{equation}
  in our formalism we need the following
  
  {\bf Lemma 5.2.1.} The following identities hold
  $$\int e^{\frac{1}{\hbar}(\frac{1}{2}\sum_{i,j=1}^na_{i,j}u_iu_j+\sum_{i=1}^nb_iu_i)}du_1...du_n=\frac{(2\pi\hbar)^{\frac{n}{2}}}{(\det A)^{\frac{1}{2}}}e^{-\frac{1}{2\hbar}\vec{b}A^{-1}\vec{b}^t}$$
  $$\vec{b}A^{-1}\vec{b}^t-c=-\frac{\det \begin{pmatrix}
A & \vec{b}^t \\
\vec{b} & c  
\end{pmatrix}}{\det A}$$
where $A=(a_{i,j})$ is a non-degenerate symmetric $n\times n$ matrix, $\vec{b}=(b_1,...,b_n)$ is a vector and $c$ is a constant.

{\bf Proof.} The first identity is a standard Gaussian integral up to multiplication by $\pm i$. The second identity is obtained by the last row and the last column expansion of the determinant in numerator in the r.h.s.
  $\square$

  {\bf Lemma 5.2.2.} The following identity holds
$$\frac{1}{2\pi\hbar}\int ue^{\frac{au^2}{\hbar}}du=\delta(a).$$

{\bf Proof.} Multiply both sides of this identity by $e^{\frac{av}{\hbar}}$ and integrate by $a$. In the r.h.s. we get 1, and in the l.h.s.
$$\frac{1}{2\pi\hbar}\int ue^{\frac{au^2+av}{\hbar}}duda=\int\delta(u^2+v)udu=1.$$
$\square$

To compute the integral (\ref{int}) we first integrate with respect to $u_1,...,u_n$ using identities from Lemma 5.2.1 and get
$$\frac{(2\pi\hbar)^{\frac{n}{2}}}{\Big(\det \Big(\sum_{i=0}^{n+1}a_{i,j,k}y_i\Big)_{1\leq j,k\leq n}\Big)^{\frac{1}{2}}}\int u_0~e^{\frac{u_0^2}{2\hbar}\frac{\det\Big(\sum_{i=0}^{n+1}a_{i,j,k}y_i\Big)_{0\leq j,k\leq n}}{\det\Big(\sum_{i=0}^{n+1}a_{i,j,k}y_i\Big)_{1\leq j,k\leq n}}}~du_0.$$

Finally, we integrate by $u_0$ using Lemma 5.2.2 and get 
$$\int u_0e^{\frac{1}{2\hbar}\sum_{i=0}^{n+1}\sum_{j,k=0}^na_{i,j,k}y_iu_ju_k}du_0...du_n=$$
$$(2\pi\hbar)^{\frac{n+2}{2}}\Big(\det \Big(\sum_{i=0}^{n+1}a_{i,j,k}y_i\Big)_{1\leq j,k\leq n}\Big)^{\frac{1}{2}}2\delta\Bigg(\det\Bigg(\sum_{i=0}^{n+1}a_{i,j,k}y_i\Bigg)_{0\leq j,k\leq n}\Bigg).$$
$\square$

One can reformulate the previous theorem using formalism from Section 3.3. Namely, consider projectively dual cones
$$C_{univ}=\{X\in Mat((n+1)\times (n+1))~\vert~X^t=X,~rk X\leq 1\}\subset\A^{\frac{(n+1)(n+2)}{2}},$$
$$\hat{C}_{univ}=\{Y\in Mat((n+1)\times (n+1))~\vert~Y^t=Y,~\det(Y)=0\}\subset(\A^{\frac{(n+1)(n+2)}{2}})^*.$$
Notice that $\dim C_{univ}=n+1$ while $\hat{C}_{univ}$ is a hypersurface. Cone $C_{univ}$ admits quadratic parametrization $X=(x_{i,j})_{0\leq i,j\leq n}$ where $x_{i,j}=u_iu_j$ (passing to the projectivization we obtain the Veronese map $\mathbb P^n\to \mathbb P^{\frac{(n+1)(n+2)}{2}-1}$ from the classical algebraic geometry). 

{\bf Theorem 5.2.2.} The Fourier transform of the density $u_0du_0...du_n$ on $C_{univ}$ is supported on $\hat{C}_{univ}$. More precisely
$$\int e^{\frac{1}{2\hbar}\sum_{j,k=0}^ny_{j,k}u_ju_k}u_0du_0...du_n=(2\pi\hbar)^{\frac{n+2}{2}}(\det(y_{j,k})_{1\leq j,k\leq n})^{\frac{1}{2}}2\delta(\det(y_{j,k})_{0\leq j,k\leq n}).$$
The proof is omitted as it is essentially the same as for Theorem 5.2.1.

One can deduce geometrically Theorem 5.2.1 from Theorem 5.2.2 in the following way. The choice of $(a_{i,j,k})$ gives a linear projection $\A^{\frac{(n+1)(n+2)}{2}}\twoheadrightarrow\A^{n+2}$. The image of $C_{univ}$ becomes a conical hypersurface in $\A^{n+2}$. The Fourier transform of the density $u_0du_0...du_n$ on the image is the restriction to the subspace $(\A^{n+2})^*\hookrightarrow (\A^{\frac{(n+1)(n+2)}{2}})^*$ of the corresponding density on $\hat{C}_{univ}$ and is supported on $(\A^{n+2})^*\cap\hat{C}_{univ}=\hat{C}$.

\ 

{\bf Remark 5.2.1.} It would be interesting to study homogeneous densities supported on conical orbit of reductive groups in finite-dimensional representations such that their Fourier transforms are supported on hypersurfaces (or cones of lower dimensions). By a similar consideration as above, this could lead to a new class of admissible conical hypersurfaces beyond those admitting quadratic parametrization.

\ 

{\bf Remark 5.2.2.} One can write system of differential equations for admissible functions $g(u_0,...,u_n)$ using results of Section 4.3, see for example (\ref{eq}). This system depends on tensor $a_{i,j,k}$. It could be holonomic or nonholonomic depending on $a_{i,j,k}$. If it is holonomic, then $rk(C)$ is finite, otherwise it is infinite.

\ 

{\bf Conjecture 5.2.1.} If $rk(C)$ is finite, then any admissible $g$ has a form
$$g(u_0,...,u_n)=\frac{Q}{P_1...P_k},~~~k=0,1,...$$
where $Q$ is a homogeneous polynomial in $u_0,...u_n$ of degree $k+1$, and $P_1,...,P_k$ are homogeneous linear polynomials in $u_0,...,u_n$.

\ 

   {\bf Remark 5.2.3.} Theorem 5.2.1 gives the lower bound $n+1$ for the rank  $rk(C)$ of  cone $C$ admitting quadratic parametrization. Experiments show that if $n\geq 3$ and tensor $a_{i,j,k}$ is generic, then this bound is attained and $rk(C)=n+1$. If $n=2$ and tensor $a_{i,j,k}$ is generic, then $rk(C)=6$, see Section 5.4. On the other hand, $rk(C)$ could be infinite for $n\geq 2$. For example, the standard smooth quadric hypersurface admits a quadratic parametrization, and its rank is infinite for $n\geq 2$. It would be interesting to find $rk(C)$ in terms of algebraic properties of tensor $a_{i,j,k}$. In particular, it would be interesting to find conditions for tensor $a_{i,j,k}$ equivalent to the property $rk(C)>n+1$.

   \ 

   {\bf Conjecture 5.2.2.} If $rk(C)$ is finite, then $rk(C)\leq (n+1)!$. 

   It would be interesting to prove this Conjecture and classify all tensors $a_{i,j,k}$ such that $rk(C)=(n+1)!$. See Section 5.4 where we construct an example of such cone for arbitrary $n$.

 \subsection{Ruled surfaces in $\P^3$}
 
 {\bf Theorem 5.3.1.} Let $\Sigma\subset\P^3$ be a surface defined parametrically by
 $$x_0=1,~~~x_1=p_1(u_2)+q_1(u_2)u_1,~~~x_2=p_2(u_2)+q_2(u_2)u_1,~~~x_3=p_3(u_2)+q_3(u_2)u_1$$
 where $p_1,p_2,p_3,q_1,q_2,q_3$ are arbitrary generic functions in one variable. Then $\Sigma$ is admissible and its rank is infinity.
 
  The corresponding admissible pairs $f(x_1,x_2),g(x_1,x_2)$ have the following parametrization
   \begin{equation}\label{rs}
   f=p_3(u_2)+q_3(u_2)u_1,~~~x_1=p_1(u_2)+q_1(u_2)u_1,~~~x_2=p_2(u_2)+q_2(u_2)u_1,
   \end{equation}
   $$g=\frac{h(u_2)}{q_1(u_2)\big(p^{\prime}_2(u_2)+q^{\prime}_2(u_2)u_1\big)-q_2(u_2)\big(p^{\prime}_1(u_2)+q^{\prime}_1(u_2)u_1\big)}$$
   where $h(u_2)$ is an arbitrary function. 
 
  {\bf Proof.}   Substituting the parametric form (\ref{rs}) in the l.h.s. of (\ref{eq0}) we get
  $$\int h(u_2)e^{\frac{1}{\hbar}(p_3(u_2)+q_3(u_2)u_1+(p_1(u_2)+q_1(u_2)u_1)y_1+(p_2(u_2)+q_2(u_2)u_1)y_2)}du_1du_2$$
  After integrating by $u_1$ we obtain
  $$2\pi\hbar\int \delta(q_3(u_2)+q_1(u_2)y_1+q_2(u_2)y_2)h(u_2)e^{\frac{1}{\hbar}(p_3(u_2)+p_1(u_2)y_1+p_2(u_2)y_2)}du_2$$
  and integrating by $u_2$ we obtain an expression of a form of the r.h.s. of (\ref{eq0}).
  $\square$

  \ 
  
  {\bf Remark 5.3.1.} Notice that the parametric representation (\ref{rs}) admits the change of variables  $u_1\to a(u_2)u_1+b(u_2),~u_2\to c(u_2)$, so the family of admissible 
  functions (\ref{rs}) is parameterized by four functions in one variable. We can set, for example, $p_1(u_2)=0, q_1(u_2)=1,q_2(u_2)=u_2$.

  \ 

  {\bf Remark 5.3.2.} For a generic ruled surface the whole space $V_f$ of admissible functions $g$ is given by (\ref{rs}). However, in special cases the space $V_f$ of admissible functions $g$ can be bigger. For example it is bigger if $f$ corresponds to a quadric. If $f(x_1,x_2)=x_1x_2$, then $V_f$ consists of functions of the form $g(x_1,x_2)=h_1(x_1)+h_2(x_2)$ where $h_1,h_2$ are arbitrary functions in one variable. It would be interesting to study in details for which ruled surfaces the space $V_f$ is bigger than those given by (\ref{rs}).

 \subsection{Steiner Roman surface}

 Hypersurfaces admitting quadratic parametrization in the case $n=2$ are equivalent to each other for the generic tensor $a_{i,j,k}$. The corresponding projective variety is called Steiner Roman surface \cite{qp}.
 
 {\bf Theorem 5.4.1.}  Let  $\Sigma\subset\P^3$ be a surface defined parametrically by 
 $$x_i=q_i(u_1,u_2),~i=0,...,3$$
 where $q_i$ are generic nonhomogeneous quadratic polynomials in $u_1,u_2$. Then $\Sigma$ is admissible and its rank is 6.

 {\bf Proof.}   It is known that any such surface is equivalent to projectivization of the affine surface defined by $x_3=f(x_1,x_2)$ where 
 $$f(x_1,x_2)=x_1x_2+\frac{x_1}{x_2}+\frac{x_2}{x_1}.$$
 In this case pair $f,g$ is admissible iff 
 $$g(x_1,x_2)=c_1+\frac{c_2}{x_1}+\frac{c_3}{x_2}+c_4~\Big(\frac{1}{x_1^2}-\frac{1}{x_2^2}\Big)+c_5~\Big(x_1-\frac{x_1}{x_2^2}\Big)+c_6~\Big(x_2-\frac{x_2}{x_1^2}\Big)$$
 where $c_1,...,c_6$ are constants. One can check this by solving the system of differential equations (\ref{eq}) with given function $f$.  $\square$

  \subsection{Generalized Steiner Roman hypersurfaces}

  Let cone $C\subset\A^{n+2}$ for $n\geq 2$ be given by
  $$x_0^{\frac{1}{2}}+x_1^{\frac{1}{2}}+...+x_{n+1}^{\frac{1}{2}}=0.$$

  {\bf Remark 5.5.1.} If $n=2$, then $C$ is equivalent to Steiner Roman surface (Section 5.4).

  {\bf Theorem 5.5.1.} The cone $C$ (and therefore its dual $\hat{C}~$) is admissible with the rank equal to $(n+1)!$.

  {\bf Proof.} The equation for the dual cone $\hat{C}$ can be written as
  $$\frac{1}{y_0}+\frac{1}{y_1}+...+\frac{1}{y_{n+1}}=0$$
  or, equivalently, as 
  $$\sum_{i=0}^{n+1}y_0...\widehat{y}_i...y_{n+1}=0.$$

  Therefore, we can choose the formal wave function in the form
  $$\tilde{g}(x_0,...,x_n)\delta\Big(x_0^{\frac{1}{2}}+x_1^{\frac{1}{2}}+...+x_{n+1}^{\frac{1}{2}}\Big)$$
  and it should be annihilated by differential operator
  $$\sum_{i=0}^{n+1}\partial_{x_0}...\widehat{\partial}_{x_i}...\partial_{x_{n+1}}.$$
  Let us make the change of variables $x_i=u_i^2,~~i=0,...,n+1$. In the new variables our wave functions reads
  $$\psi=g(u_0,...,u_n)\delta(u_0+u_1+...+u_{n+1})$$
  and the differential operator (up to a coefficient and after multiplication by $u_0u_1...u_{n+1}$) takes the form
  $$D=\sum_{i=0}^{n+1}u_i\partial_{u_0}...\widehat{\partial}_{u_i}...\partial_{u_{n+1}}.$$
  The equation $D(\psi)=0$ is equivalent to the following system of differential equations for the function $g(u_0,...,u_n)$
  $$\sum_{i=0}^{n}u_i\partial_{u_i}g+(n+1)g=0$$
  $$\sum_{0\leq i_1 < i_2\leq n}(u_{i_1}+u_{i_2})\partial_{u_{i_1}}\partial_{u_{i_2}}g+n\sum_{i=0}^n\partial_{u_i}g=0$$
  $$........................................................................................................................................$$
  \begin{equation}\label{hols}
  \sum_{0\leq i_1<...<i_k\leq n}(u_{i_1}+...+u_{i_k})\partial_{u_{i_1}}...\partial_{u_{i_k}}g+(n+2-k)\sum_{0\leq j_1<...<j_{k-1}\leq n}\partial_{u_{j_1}}...\partial_{u_{j_{k-1}}}g=0
  \end{equation}
  $$........................................................................................................................................$$
  $$(u_0+...+u_n)\partial_{u_0}...\partial_{u_n}g+\sum_{j=0}^n\partial_{u_0}...\widehat{\partial}_{u_j}...\partial_{u_n}g=0.$$
  To finish the proof we need the following

  {\bf Lemma 5.5.1.} The system (\ref{hols}) is holonomic and has no more than $(n+1)!$-dimensional space of solutions at generic point $(u_0,...,u_n)$.

  {\bf Proof.} The system (\ref{hols}) defines a cyclic $D$-module $M$ generated by $g$. The filtration $D^{\leq i}$ by the degree of differential operators, after applying to the generator $g$ induces a good filtration on $M$. The associated graded space $grM$ is a graded module over polynomial ring ${\bf k}[u_0,...,u_n,v_0,...,v_n]$ where $\deg u_i=0,~\deg v_i=1$ and $v_i$ are images of $\partial_{u_i}$ in the polynomial ring. The module $grM$ is also cyclic and defined by equations
  $$\Big(u_0v_0+...+u_nv_n\Big)h=0$$
  $$\sum_{0\leq i_1 < i_2\leq n}(u_{i_1}+u_{i_2})v_{i_1}v_{i_2}h=0$$
  $$.........................................................................$$
  $$(u_0+...+u_n)v_0...v_nh=0.$$
  One can check that for generic $u_i$ (e.g. $u_0=...=u_n=1$) the above system gives a 0-dimensional scheme which is a complete intersection of length $(n+1)!$, the product of degrees in $v_i$ of our equations.
  $\square$

  {\bf Lemma 5.5.2.} Let $\{0,1,...,n+1\}=P_1\sqcup...\sqcup P_{m+1},~~m=1,...,n+1$ be an arbitrary partition. Set $u_{n+1}=-u_0-...-u_n$.
  Then the following functions 
  $$g=\frac{1}{u_0...u_{n+1}}\cdot \frac{u_{i_1}...u_{i_{m+1}}}{\sum_{\alpha\in P_1}u_{\alpha}\cdot \sum_{\alpha\in P_2}u_{\alpha}\cdot ...\cdot \sum_{\alpha\in P_m}u_{\alpha}},~~~i_1\in P_1,...,i_{m+1}\in P_{m+1}$$
  satisfy the system (\ref{hols}). The vector space spent by these functions is $(n+1)!$-dimensional.

  {\bf Proof.} This can be proved by induction as follows. Without loss of generality we can set $P_1=\{0,...,k\}$ and $i_1=k$ because of $S_{n+1}$ action on the system (\ref{hols}). Write 
  $$g(u_0,...,u_n)=\frac{g_1(u_{k+1},...,u_n)}{u_0...u_{k-1}(u_0+...+u_k)}$$
  and check that $g$ satisfies (\ref{hols}) iff $g_1$ satisfies similar system with $n-k$ variables $u_{k+1},...,u_n$.
  $\square$

  {\bf Example 5.5.1.} Let $n=4$, $P_1=\{0,1\},~P_2=\{2,3\},~P_3=\{4,5\}$ and choose $0\in P_1,~2\in P_2,~4\in P_3$. We have
  $$g=\frac{1}{u_0u_1u_2u_3u_4u_5}\cdot \frac{u_0u_2u_4}{(u_0+u_1)(u_2+u_3)}.$$

  {\bf Remark 5.5.2.} We have $\deg C=2^n$ and $\deg\hat{C}=n+1$. In particular, $\hat{C}$ is not equivalent to $C$ and gives another example of admissible cone with the rank equal to $(n+1)!$.

 \subsection{Kummer surfaces in $\P^3$}

 Recall that the Kummer surface is a quartic surface in $\P^3$ with 16 double points. Abstractly it is the quotient of a principally polarized two dimensional abelian variety by the antipodal involution. 

 {\bf Theorem 5.6.1.} Any Kummer surface is admissible of rank 6.

 {\bf Proof.} We need to check (\ref{eqalg}) for the Kummer surfaces. Choose equation for Kummer surface in the form \cite{ks}
 $$F_1=(x_0^2+x_1^2+x_2^2+x_3^2+a(x_0x_1+x_2x_3)+b(x_0x_2+x_1x_3)+c(x_0x_3+x_1x_2))^2$$
 $$-4(a^2+b^2+c^2-abc-4)x_0x_1x_2x_3.$$
 Computing equation for dual polynomial we get
 $$F_2=(4a-2bc)(\partial_{x_0}\partial_{x_1}\partial_{x_2}^2+\partial_{x_0}^2\partial_{x_2}\partial_{x_3}+\partial_{x_1}^2\partial_{x_2}\partial_{x_3}+\partial_{x_0}\partial_{x_1}\partial_{x_3}^2)+$$
 $$(4b-2ac)(\partial_{x_0}\partial_{x_1}^2\partial_{x_2}+\partial_{x_0}^2\partial_{x_1}\partial_{x_3}+\partial_{x_1}\partial_{x_2}^2\partial_{x_3}+\partial_{x_0}\partial_{x_2}\partial_{x_3}^2)+$$
 $$(4c-2ab)(\partial_{x_0}^2\partial_{x_1}\partial_{x_2}+\partial_{x_0}\partial_{x_1}^2\partial_{x_3}+\partial_{x_0}\partial_{x_2}^2\partial_{x_3}+\partial_{x_1}\partial_{x_2}\partial_{x_3}^2)+$$
 $$(a^2-4)(\partial_{x_0}^2\partial_{x_1}^2+\partial_{x_2}^2\partial_{x_3}^2)+(b^2-4)(\partial_{x_0}^2\partial_{x_2}^2+\partial_{x_1}^2\partial_{x_3}^2)+(c^2-4)(\partial_{x_0}^2\partial_{x_3}^2+\partial_{x_1}^2\partial_{x_2}^2)+$$
 $$(4abc-2a^2-2b^2-2c^2-8)\partial_{x_0}\partial_{x_1}\partial_{x_2}\partial_{x_3}.$$
 Let us choose $g$ in the form
 $$g=x_1^{\frac{1}{2}}x_2^{\frac{1}{2}}\Big((b\sqrt{c^2-4}+2a-bc)(x_0^2+cx_0x_3+x_3^2)+2\sqrt{c^2-4}(x_0x_2+x_1x_3)+$$
 $$+(c\sqrt{c^2-4}-c^2+4)(x_0x_1+x_2x_3)\Big)^{\frac{1}{2}}.$$
 One can check by direct calculation that
 $$F_2 \Big( g\cdot\delta(F_1)\Big)=0$$
 where $F_2$ is a differential operator and $g,F_1$ are expression written above. Notice that $F_1,F_2$ are invariant with respect to the  action of the group of permutation $S_4$ acting on $x_0,...,x_3$ and simultaneously on $a,b,c$. Applying this action to $g$ we obtain other admissible functions. We can change sign of the square root $\sqrt{c^2-4}$ as well. One can check that in this way we obtain a six-dimensional space of admissible functions $g$. On the other hand, we know (see Section 7) that if an admissible surface is not ruled, then $\dim V_f\leq 6$.
 $\square$

 The space $V_f$ for a Kummer surface can be constructed as follows. Let $\mathcal A_2$ be a two-dimensional principally polarized abelian variety. Let $\sigma: \mathcal A\to\mathcal A$ be its involution with a fixed point. Construct the corresponding  Kummer surface by $\mathcal K =\mathcal A_2 /\sigma.$ Let $\mathcal L$ be a line bundle of degree 2 on $\mathcal A_2$ (hence with 4 sections) and such that  $\mathcal L^{\sigma}=\mathcal L$. It is known that $\sigma$ acts as identity operator on the space $\Gamma(\mathcal L)$ of sections of $\mathcal L$. We can identify $x_0,...,x_3$ with a basis in $\Gamma(\mathcal L)$. It is also known that $\dim\Gamma(\mathcal L^2)=16$ and $\Gamma(\mathcal L^2)=\Gamma^+\oplus\Gamma^-$ where $\sigma$ acts as identity on $\Gamma^+$ and as negative identity on $\Gamma^-$. We have $\dim\Gamma^+=10$ and $x_ix_j,~0\leq i\leq j\leq 3$ is a basis in $\Gamma^+$. The space $V_f$ can be identified with $\Gamma^-$. 

 These spaces of sections can be constructed explicitly over $\C$ using theta functions in two variables.

 Define the functions $\theta_{a,b}(u_1,u_2),~a,b\in\Z/2$ by
 $$\theta_{00}(u_1,u_2)=\sum_{i,j\in\Z}e^{2\pi\sqrt{-1}~\big(2iu_1+2ju_2+i(i-1)\tau_{11}+2ij\tau_{12}+j(j-1)\tau_{22}\big)},$$
 $$\theta_{10}(u_1,u_2)=\sum_{i,j\in\Z}e^{2\pi\sqrt{-1}~\big((2i+1)u_1+2ju_2+i^2\tau_{11}+(2i+1)j\tau_{12}+j(j-1)\tau_{22}\big)},$$
 $$\theta_{01}(u_1,u_2)=\sum_{i,j\in\Z}e^{2\pi\sqrt{-1}~\big(2iu_1+(2j+1)u_2+i(i-1)\tau_{11}+i(2j+1)\tau_{12}+j^2\tau_{22}\big)},$$
 $$\theta_{11}(u_1,u_2)=\sum_{i,j\in\Z}e^{2\pi\sqrt{-1}~\big((2i+1)u_1+(2j+1)u_2+i^2\tau_{11}+(2ij+i+j)\tau_{12}+j^2\tau_{22}\big)}$$
 where $\tau_{00},\tau_{12},\tau_{22}$ are constants such that the matrix $\begin{pmatrix} \tau_{11}&\tau_{12}\\ \tau_{11}&\tau_{22}\end{pmatrix}$  
 belongs to Siegel upper half-space of genus two, i.e. its imaginary part is positive. The functions $\theta_{a,b}(u_1,u_2)$ are holomorphic and satisfy the following periodicity 
 and quasiperiodicity  properties
 $$\theta_{a,b}(u_1+1,u_2)=\theta_{a,b}(u_1,u_2),~~~\theta_{a,b}(u_1,u_2+1)=\theta_{a,b}(u_1,u_2),$$
 $$\theta_{a,b}(u_1+\tau_{11},u_2+\tau_{12})=e^{-2\pi\sqrt{-1}~2u_1}\theta_{a,b}(u_1,u_2),~~~\theta_{a,b}(u_1+\tau_{12},u_2+\tau_{22})=e^{-2\pi\sqrt{-1}~2u_2}\theta_{a,b}(u_1,u_2).$$
 Moreover, any holomorphic function with these properties is a linear combination of $\theta_{a,b}(u_1,u_2),$ $a,b\in\Z/2$.
 
 We can define Kummer surface and corresponding functions $f(x_1,x_2),g(x_1,x_2)$ parametrically by 
 $$f=\frac{\theta_{00}(u_1,u_2)}{\theta_{11}(u_1,u_2)},~x_1=\frac{\theta_{10}(u_1,u_2)}{\theta_{11}(u_1,u_2)},~x_2=\frac{\theta_{01}(u_1,u_2)}{\theta_{11}(u_1,u_2)},$$
 $$g=\frac{1}{\Delta}~\Bigg(c_1\Bigg(\frac{\theta_{00}(u_1,u_2)}{\theta_{11}(u_1,u_2)}\Bigg)_{u_1}+c_2\Bigg(\frac{\theta_{10}(u_1,u_2)}{\theta_{11}(u_1,u_2)}\Bigg)_{u_1}+
 c_3\Bigg(\frac{\theta_{01}(u_1,u_2)}{\theta_{11}(u_1,u_2)}\Bigg)_{u_1}+$$
 $$c_4\Bigg(\frac{\theta_{00}(u_1,u_2)}{\theta_{11}(u_1,u_2)}\Bigg)_{u_2}+c_5\Bigg(\frac{\theta_{10}(u_1,u_2)}{\theta_{11}(u_1,u_2)}\Bigg)_{u_2}+
 c_6\Bigg(\frac{\theta_{01}(u_1,u_2)}{\theta_{11}(u_1,u_2)}\Bigg)_{u_2}\Bigg)$$
 where $c_1,...,c_6$ are constants, indexes $u_1,u_2$ stand for partial derivatives and
 $$\Delta=\Bigg(\frac{\theta_{10}(u_1,u_2)}{\theta_{11}(u_1,u_2)}\Bigg)_{u_1}~\Bigg(\frac{\theta_{01}(u_1,u_2)}{\theta_{11}(u_1,u_2)}\Bigg)_{u_2}-
 \Bigg(\frac{\theta_{01}(u_1,u_2)}{\theta_{11}(u_1,u_2)}\Bigg)_{u_1}~\Bigg(\frac{\theta_{10}(u_1,u_2)}{\theta_{11}(u_1,u_2)}\Bigg)_{u_2}.$$
 Then pair $f,g$ is admissible.

 {\bf Remark 5.6.1.} Kummer surfaces have many degenerations which are not Kummer surfaces anymore but also give examples of admissible surfaces. For example, Steiner Roman surface can be obtained in this way.

 {\bf Remark 5.6.2.} It would be interesting to find a more transparent proof of Theorem 5.6.1, for example, using parametrization of Kummer surface by theta functions or using relation of Kummer surface with the quadratic line complex. Notice that our space of sections $\Gamma^-$ described above has a basis $y_1,...,y_6$ such that
 $$y_1^2+...+y_6^2=0,$$
  $$\lambda_1y_1^2+...+\lambda_6y_6^2=0,$$
   $$\lambda_1^2y_1^2+...+\lambda_6^2y_6^2=0$$
   where $\lambda_1,...,\lambda_6\in\C$ are constants. Therefore, the space $\Gamma^-$ gives an embedding of the quadratic line complex to $\P^5$.
 
 \subsection{Extensions of admissible pairs and families of surfaces of degree four}

 {\bf Theorem 5.7.1.} Let a pair of functions $\tilde{f}(x_1,...,x_k,y_{k+1},...,y_n),~\tilde{g}(x_1,...,x_k,y_{k+1},...,y_n)$ be admissible as functions in variables $x_1,...,x_k$ where $y_{k+1},...,y_n$ are considered as parameters, and let this pair be also admissible as functions in variables $y_{k+1},...,y_n$ where $x_1,...,x_k$ are considered as parameters. Let us construct a new pair of functions $f(x_1,...,x_n),~g(x_1,...,x_n)$ by doing the Fourier transform with respect to variables $y_{k+1},...,y_n$ and using our admissibility assumption
 $$\int \tilde{g}(x_1,...,x_k,y_{k+1},...,y_n)e^{\frac{\tilde{f}(x_1,...,x_k,y_{k+1},...,y_n)-x_{k+1}y_{k+1}-...-x_ny_n}{\hbar}}dy_{k+1}...dy_n=$$
 $$(2\pi \hbar)^{\frac{n-k}{2}} g(x_1,...,x_n)e^{\frac{f(x_1,...,x_n)}{\hbar}}.$$
 Then the pair of functions in $n$ variables  $f,g$ is admissible.

 {\bf Proof.} Notice that the pair $f,g$ is also admissible with respect to variables $x_{k+1},...,x_n$ as the Fourier transform of an admissible pair. We compute the integral in the l.h.s. of (\ref{eq0}) in two steps so that in each step the integral is one loop exact. First, we integrate with respect to $x_{k+1},...,x_n$ using the admissibility with respect to these variables. Second, we integrate with respect to $x_1,...,x_k$ using admissibility of $\tilde{f},\tilde{g}$ with respect to these variables. $\square$

 {\bf Example 5.7.1.} Let us define functions $\tilde{f}(x_1,y_2),~\tilde{g}(x_1,y_2)$ by
 $$\tilde{f}(x_1,y_2)=\sqrt{x_1^2-1}\sqrt{y_2^2-1}+(a_1x_1+a_2)\sqrt{y_2^2-1}+(a_3y_2+a_4)\sqrt{x_1^2-1}+a_5x_1y_2,$$
 $$\tilde{g}(x_1,y_2)=\frac{c_1}{\sqrt{x_1+1}\sqrt{y_2+1}}+\frac{c_2}{\sqrt{x_1-1}\sqrt{y_2+1}}+\frac{c_3}{\sqrt{x_1+1}\sqrt{y_2-1}}+\frac{c_4}{\sqrt{x_1-1}\sqrt{y_2-1}}$$
 where $a_1,...,a_5,c_1,...,c_4$ are constants. 
 It follows from results of Section 6 (see Theorem 6.1) that the pair $\tilde{f},\tilde{g}$ is admissible with respect to variable $x_1$ and with respect to variable $y_2$ (but it is not admissible as a pair of functions in two variables). After the Fourier transform with respect to $y_2$ we get the following admissible pair of functions in two variables
 $$f(x_1,x_2)=a_4\sqrt{x_1^2-1}+$$
 $$\sqrt{x_2-(a_1+a_5)x_1-a_2-(a_3+1)\sqrt{x_1^2-1}}\sqrt{x_2+(a_1-a_5)x_1+a_2-(a_3-1)\sqrt{x_1^2-1}},$$
 $$g(x_1,x_2)=\frac{\frac{c_1}{\sqrt{x_1+1}}+\frac{c_2}{\sqrt{x_1-1}}}{\sqrt{x_2-(a_1+a_5)x_1-a_2-(a_3+1)\sqrt{x_1^2-1}}}+$$
 $$\frac{\frac{c_3}{\sqrt{x_1+1}}+\frac{c_4}{\sqrt{x_1-1}}}{\sqrt{x_2+(a_1-a_5)x_1+a_2-(a_3-1)\sqrt{x_1^2-1}}}.$$

 A similar example of an admissible pair coming from Theorem 5.7.1 is given below, where explicit formulas for $\tilde{f},\tilde{g}$ are omitted. 
 
 {\bf Theorem 5.7.2.} Define functions $f(x_1,x_2),g(x_1,x_2)$ parametrically as 
 $$f=u_1u_2^4 + 2r_1r_2u_1u_2^3 + (r_1^2r_2^2 + 2)u_1u_2^2 + 2u_1u_2r_1r_2 + u_1 + \frac{1}{u_1},~x_1=u_2^2r_2 - u_2r_1 + r_2 + \frac{1}{u_1}, ~x_2=u_2^2$$
 $$g=\sqrt{u_1}~C_1 + \sqrt{u_1}~\frac{C_2}{u_2} + u_1^{3/2}(r_1r_2u_2 + u_2^2 + 1)C_3 + u_1^{3/2}(r_1r_2u_2 + u_2^2 + 1)\frac{C_4}{u_2}$$
 where $r_1,r_2,C_1,C_2,C_3,C_4$ are constants. Then the pair $f(x_1,x_2),g(x_1,x_2)$ is admissible and its rank is 4.
 
 {\bf Proof.} Since $f$ is an algebraic function, we just need to check (\ref{eq}) which can be done by a direct calculation. Alternatively, notice that for any fixed $x_2$ the pair $f,g$ is admissible as pair of functions in one variable $x_1$.  Moreover, the integral 
 $\int g(x_1,x_2) e^{\frac{f(x_1,x_2)+x_1y_1+x_2y_2}{\hbar}}dx_1dx_2$ can be computed in two steps: first with respect to $x_1$, and than with respect to $x_2$. In each step the 
 integral with respect to one variable is one loop exact. $\square$

 {\bf Remark 5.7.1.} One can check that Example 5.7.1 and Theorem 5.7.2 both give a two parametric family of admissible cones in $\P^3$ of degree 4. In both cases these cones have non-isolated singularities. We have not checked if these cones are projectively equivalent.

 {\bf Remark 5.7.2.} One can generalize the construction from Theorem 5.7.1 considering for example functions in three or more groups of variables. It might be interesting to study admissible pairs obtained in this way. On the other hand, these admissible pairs look degenerate in a certain sense.

 \subsection{Toric hypersurfaces}

 Admissible pairs of the form $f=x_1^{a_1}...x_n^{a_n},~g=x_1^{b_1}...x_n^{b_n}$ were completely described in \cite{efp}. Here we just reformulate their results in terms of our formalism.
 
 {\bf Theorem 5.8.1.} Let $C\subset\P^{n+1}$ be a toric projective hypersurface defined by 
 \begin{equation}\label{tor}
 x_0^{p_0}...x_{m-1}^{p_{m-1}}=x_{m}^{p_{m}}...x_{n+1}^{p_{n+1}}
 \end{equation}
 where $p_0,...,p_{n+1}>0$ are integers and $p_0+...+p_{m-1}=p_{m+}+...+p_{n+1}$.  Then $C$ is admissible if there exist numbers $r_0,...,r_{n+1}\in\Q$ such the following identity holds as an identity for functions in one variable $t$
 \begin{equation}\label{gam}
 \Gamma(p_0t+r_0)...\Gamma(p_{m-1}t+r_{m-1})=ae^{bt}\Gamma(p_{m}t+r_{m})...\Gamma(p_{n+1}t+r_{n+1})
 \end{equation}
 where $a,b$ are constants. 
 
  We expect that the rank $C$ is equal to the number of vectors $(r_0,...,r_{n+1})\in\Q^{n+2}$ up to translation of the form $r_i\to r_i+p_iu,~i=0,...,n+1$ for some $u$, such that the identity (\ref{gam})
   holds.

   {\bf Remark 5.8.1.} The transcendental condition (\ref{gam}) is equivalent to a combinatorial one: the union of sets of poles (with multiplicities) in the l.h.s. and in the r.h.s. of equation (\ref{gam}) coincide:

   $$\bigcup_{i=0}^{m-1} \frac{1}{p_i}(\Z_{\geq 0}+r_i)=\bigcup_{i=m}^{n+1} \frac{1}{p_i}(\Z_{\geq 0}+r_i).$$

   Notice that here we multiply $t$ by $-1$.
  
  {\bf Proof.}  Let $f(x_1,...,x_n)=x_1^{a_1}...x_n^{a_n}$ where $a_1,...,a_n\in\Q\setminus 0$. Let us set 
  $$g(x_1,...,x_n)=x_1^{b_1}...x_n^{b_n}.$$
  Let
   \begin{equation}\label{I}
   I(\vec{y})=\int x_1^{b_1}...x_n^{b_n}e^{\frac{1}{\hbar}(x_1^{a_1}...x_n^{a_n}+x_1y_1+...+x_ny_n)}.
   \end{equation}

   {\bf Lemma 5.8.1.} The following equality of the formal wave functions holds
   \begin{equation}\label{exp}
   e^{\frac{x}{\hbar}}=\int \frac{(\frac{x}{\hbar})^{\frac{t}{\hbar}}}{\Gamma(\frac{t}{\hbar}+1)}\frac{dt}{\hbar}.
   \end{equation}

   {\bf Proof.} First, let us show that the r.h.s. of this equation makes sense as a formal wave function. Let 
   $$S(t)=\exp\Big(-\frac{\zeta(-1)}{1}t-\frac{\zeta(-3)}{3}t^3-\frac{\zeta(-5)}{5}t^5-...\Big)\in\Q[[t]]$$
   where $\zeta$ stands for the Riemann zeta function. Notice that 
   \begin{equation}\label{S}
   S(t)S(-t)=1.
   \end{equation}
   Recall Stirling's formula
   $$\Gamma(t+1)\sim (2\pi t)^{\frac{1}{2}}t^te^{-t}S(t^{-1})$$
   where $\sim$ means the asymptotic expansion for $t\to\infty$. Since the r.h.s. of (\ref{exp}) is also understood as an asymptotic expansion, we can replace the gamma function with the r.h.s. of Stirling's formula and rewrite the equation (\ref{exp}) as 
   $$(2\pi\hbar)^{\frac{1}{2}} e^{\frac{x}{\hbar}}=\int t^{-\frac{1}{2}}S\Big(-\frac{\hbar}{t}\Big)e^{\frac{t}{\hbar}(1+\ln x-\ln t)}dt.$$
   The r.h.s. of this integral can be computed by the  expansion at the critical point $t=x$.

   Now let us prove the identity (\ref{exp}). Let $\phi_l(x)$ (resp. $\phi_r(x)$) be the l.h.s. (resp. the r.h.s.) of (\ref{exp}). One can check that these functions satisfy the same differential equations $\partial_x\phi_l(x)=\frac{x}{\hbar}\phi_l(x),~\partial_x\phi_r(x)=\frac{x}{\hbar}\phi_r(x)$. Therefore, we have $\phi_l(x)=c\phi_r(x)$ where $c$ is independent of $x$. To compute $c$ we just compute the first term of (\ref{exp}) by the stationary phase method. $\square$

   Using Lemma 5.8.1 we replace $e^{\frac{1}{\hbar}x_1^{a_1}...x_n^{a_n}}$ in $I(\vec{y})$ with the corresponding r.h.s. of the identity (\ref{exp}) and obtain
    \begin{equation}\label{I1}
   I(\vec{y})=\int\prod_{i=1}^nx_i^{b_i+\frac{a_it}{\hbar}}e^{\frac{x_iy_i}{\hbar}}~\frac{1}{\Gamma(\frac{t}{\hbar}+1)\hbar^{\frac{t}{\hbar}}}~\frac{dt}{\hbar}~dx_1...dx_n.
   \end{equation}

   {\bf Lemma 5.8.2.} Let $a>0$. The following equalities between formal wave functions hold
   $$\int x^{b+\frac{a}{\hbar}}e^{\frac{xy}{\hbar}}dx=\Big(-\frac{y}{\hbar}\Big)^{-b-\frac{a}{\hbar}-1}\Gamma(b+\frac{a}{\hbar}+1),$$
    $$\int x^{b-\frac{a}{\hbar}}e^{\frac{xy}{\hbar}}dx=\Big(\frac{y}{\hbar}\Big)^{-b+\frac{a}{\hbar}-1}\frac{2\pi i}{\Gamma(-b+\frac{a}{\hbar})}.$$

    {\bf Proof.} For the first equation notice that if $x>0,~y<0,~\hbar>0$, then we have the usual Euler integral representation of the gamma function in the form 
    $$\int_0^{\infty} x^{b+\frac{a}{\hbar}}e^{\frac{xy}{\hbar}}dx=\Big(-\frac{y}{\hbar}\Big)^{-b-\frac{a}{\hbar}-1}\Gamma(b+\frac{a}{\hbar}+1)$$
    which implies our equality for the formal wave functions. The second equation makes sense only on the level of wave functions and follows from the first one and the equation (\ref{S}).
  $\square$
  
  Using these formulas we rewrite $I(\vec{y})$ in (\ref{I1}) as
  $$I(\vec{y})=\int\frac{1}{\Gamma(\frac{t}{\hbar}+1)\hbar^{\frac{t}{\hbar}}}\prod_{i=1}^n\Big(\pm\frac{y_i}{\hbar}\Big)^{-b_i-\frac{a_it}{\hbar}-1}\cdot \frac{\prod_{i=1}^k\Gamma(b_i+\frac{a_it}{\hbar}+1)}{\prod_{i=k+1}^n\Gamma(-b_i-\frac{a_it}{\hbar})}\cdot \frac{dt}{\hbar}$$
  where we assume $a_1,...,a_k>0,~a_{k+1},...,a_n<0$. It follows from this expression that the integral $I(\vec{y})$ is one loop exact if there exist $a,b\in\Q$, $a\ne 0$ such that 
  $$\frac{1}{\Gamma(\frac{t}{\hbar}+1)}\cdot\frac{\prod_{i=1}^k\Gamma(b_i+\frac{a_it}{\hbar}+1)}{\prod_{i=k+1}^n\Gamma(-b_i-\frac{a_it}{\hbar})}=c_1c_2^t\Gamma(b+\frac{at}{\hbar}+1)\text{~~~or~~~}\frac{c_1c_2^t}{\Gamma(-b-\frac{at}{\hbar})}.$$
  Notice that we must have $a=a_1+...+a_n-1$. Reformulating this condition in terms of toric hypersurface (\ref{tor}) we obtain the condition (\ref{gam}).
   $\square$

   The statement and the proof of this theorem can be generalized to the field $\R$ instead of $\Q$.
  
  {\bf Remark 5.8.2.} If an identity of the form (\ref{gam}) holds, then it holds by virtue of  Gauss multiplication formula. This can be shown by considering poles of the l.h.s. and the r.h.s.
   of (\ref{gam}).
   
   {\bf Example 5.8.1.}  Let $C\subset\P^{3}$ be a toric projective surface defined by 
   $$x_0x_1x_2=x_3^3.$$
   Any vector $(r_0,r_1,r_2,r_3)\in\R^4$ such that
   $$\Gamma(t+r_0)\Gamma(t+r_1)\Gamma(t+r_2)=ae^{bt}\Gamma(3t+r_3)$$
   is equal to $(0,\frac{1}{3},\frac{2}{3},0)$ up to translation by $(u,u,u,3u)$ and permutation of $r_0,r_1,r_2$. Therefore, $C$ is admissible with  rank $3!=6$.

   {\bf Example 5.8.2.} More generally, let $C\subset\P^{n+1}$ be a toric projective surface defined by 
   $$x_0...x_n=x_{n+1}^{n+1}.$$
   Then its rank is equal to $(n+1)!$.
   
    {\bf Example 5.8.3.}  Let $C\subset\P^{3}$ be a toric projective surface defined by 
   $$x_0x_1x_2^2=x_3^4.$$
   Any vector $(r_0,r_1,r_2,r_3)\in\R^4$ such that
   $$\Gamma(t+r_0)\Gamma(t+r_1)\Gamma(2t+r_2)=ae^{bt}\Gamma(4t+r_3)$$
   is equal to either $(0,\frac{1}{2},\frac{1}{2},0)$ or $(\frac{1}{4},\frac{3}{4},0,0)$ up to translation by $(u,u,2u,4u)$ and permutation of $r_0,r_1$. Therefore, $C$ is admissible with rank $2\cdot 2=4$.
   
    {\bf Example 5.8.4.}  Let $C\subset\P^{3}$ be a toric projective surface defined by 
   $$x_0x_1^p=x_2x_3^p$$
   for some $p>0$.
   Any vector $(r_0,r_1,r_2,r_3)\in\R^4$ such that
   $$\Gamma(t+r_0)\Gamma(pt+r_1)=ae^{bt}\Gamma(t+r_2)\Gamma(pt+r_3)$$
   is equal to $(0,v,0,v)$ for arbitrary $v$,  up to translation by $(u,pu,u,pu)$. Therefore, $C$ is admissible with infinite rank.
   
   {\bf Remark 5.8.3.} One can show that if $n=2$, then any admissible toric surface of the form (\ref{tor}) is given by examples 5.8.1, 5.8.3, 5.8.4 (up to permutations of $x_0,...,x_3$ and 
   multiplication of $p_0,...,p_3$ by a common factor.

   \subsection{Segre cubic in $\P^4$}

   {\bf Theorem 5.9.1.} Define a cone in $\A^5$ by 
   $$x_0x_1x_2+x_1x_2x_3+x_2x_3x_4+x_3x_4x_0+x_4x_0x_1=0.$$
   This cone is admissible with rank 24.

   {\bf Proof.} Let 
   $$F_1=x_0x_1x_2+x_1x_2x_3+x_2x_3x_4+x_3x_4x_0+x_4x_0x_1.$$
   Computing dual polynomial we obtain
   $$F_2=\partial_{x_0}^2\partial_{x_1}^2+\partial_{x_1}^2\partial_{x_2}^2+\partial_{x_2}^2\partial_{x_3}^2+\partial_{x_3}^2\partial_{x_4}^2+\partial_{x_4}^2\partial_{x_0}^2+$$
   $$+2\partial_{x_0}\partial_{x_1}\partial_{x_2}\partial_{x_3}+2\partial_{x_1}\partial_{x_2}\partial_{x_3}\partial_{x_4}+2\partial_{x_2}\partial_{x_3}\partial_{x_4}\partial_{x_0}+2\partial_{x_3}\partial_{x_4}\partial_{x_0}\partial_{x_1}+2\partial_{x_4}\partial_{x_0}\partial_{x_1}\partial_{x_2}$$
   $$-2\partial_{x_0}\partial_{x_1}^2\partial_{x_2}-2\partial_{x_1}\partial_{x_2}^2\partial_{x_3}-2\partial_{x_2}\partial_{x_3}^2\partial_{x_4}-2\partial_{x_3}\partial_{x_4}^2\partial_{x_0}-2\partial_{x_4}\partial_{x_0}^2\partial_{x_1}.$$
   Let $g$ be a cyclic permutation of one of the following expressions
   $$\sqrt{\frac{x_0x_1x_2}{x_3x_4}},~\sqrt{\frac{x_0(x_1+x_3)(x_1+x_4)}{x_1x_2}},~\sqrt{\frac{x_3(x_2+x_0)(x_2+x_4)}{x_1x_2}},$$
   $$\sqrt{\frac{x_0x_4(x_0+x_2)}{x_1(x_0+x_3)}},~\sqrt{\frac{x_0x_1(x_0+x_3)}{x_4(x_0+x_2)}},~\sqrt{\frac{x_1x_2x_3}{(x_2+x_4)(x_2+x_0)}}.$$
   Thus we have 30 possible functions $g$. One can check that the vector space spent by these functions is 24-dimensional and 
   $$F_2\Big(g\delta(F_1)\Big)=0$$
   for any such $g$. One can also check that the system for $g$ is holonomic and has a 24-dimensional space of solutions at a generic point.
   $\square$

   {\bf Remark 5.9.1.} Let $H=\det ((\partial_{x_i}\partial_{x_j}F_1)_{0\leq i,j\leq4})$. One can check that divisor on Segre surface $F_1=0$ defined by $H=0$ is equivalent to the divisor defined by 
   $$x_0x_1x_2x_3x_4(x_0+x_2)(x_1+x_3)(x_2+x_4)(x_3+x_0)(x_4+x_1)=0.$$

 \section{Classification of admissible pairs of functions in one variable}
 
 In the case of functions $f(x),g(x)$ in one variable, the classification of admissible pairs $f,g$ can be done directly, by equating the higher loop contributions to zero and solving the 
 corresponding system of differential equations. 
 
 {\bf Lemma 6.1.} If $f(x)$ is admissible, then it has one of the following forms:
 
 {\bf 1)}  $f(x)=(a_0+a_1x+a_2x^2)^{\frac{1}{2}}+b_0+b_1x$ where $a_0,a_1,a_2,b_0,b_1$ are arbitrary constants such that the polynomial $a_0+a_1x+a_2x^2$ has distinct 
 roots. In particular, if $a_2=0$, then $a_1\ne 0$ (root at infinity has multiplicity at most one).
 
 {\bf 2)} $f(x)=\frac{1}{a_0+a_1x}+b_0+b_1x$ where $a_0,a_1,b_0,b_1$ are arbitrary constants and $a_1\ne 0$.

  {\bf 3)} $f(x)=b_0+b_1x+b_2x^2$ where $b_0,b_1,b_2$ are arbitrary constants and $b_2\ne 0$.
 
 {\bf Proof.} Equating to zero the terms $A_2,A_3,A_4$ of the expansion (\ref{per}) one obtains a system of polynomial differential equations for $f(x),g(x)$.  Assuming 
 that $f^{\prime\prime}(x),g(x)\ne 0$ one can check by direct computations that this system is equivalent to the following:
 \begin{equation} \label{sys1}
 9f^{\prime\prime}(x)^2f^{(5)}(x)-45f^{\prime\prime}(x)f^{\prime\prime\prime}(x)f^{(4)}(x)+40f^{\prime\prime\prime}(x)^3=0,
 \end{equation} 
 $$12f^{\prime\prime}(x)^2\cdot g^{\prime\prime}(x)-12f^{\prime\prime}(x)f^{\prime\prime\prime}(x)\cdot g^{\prime}(x)+(5f^{\prime\prime\prime}(x)^2-3f^{\prime\prime}(x)f^{(4)}(x))\cdot g(x)=0.$$
 The first equation can be written as 
 $$\frac{d^3}{dx^3}\Big(f^{\prime\prime}(x)^{-\frac{2}{3}}\Big)=0$$
 which gives 
 $$f^{\prime\prime}(x)=(a_0+a_1x+a_2x^2)^{-\frac{3}{2}}$$
 where $a_0,a_1,a_2$ are arbitrary constants such that the vector $(a_0,a_1,a_2)$ is non-zero.  After the integration we obtain the statement of Lemma. $\square$
 
 {\bf Theorem 6.1.} The following pairs $f(x),g(x)$ are admissible:
 
 {\bf 1)} $f(x)=\sqrt{x^2-1},~~~g(x)=\frac{c_1}{\sqrt{x+1}}+\frac{c_2}{\sqrt{x-1}}$.
 
 {\bf 2)} $f(x)=\frac{1}{x},~~~g(x)=c_1x^{-\frac{1}{2}}+c_2x^{-\frac{3}{2}}$.
 
 {\bf 3)} $f(x)=x^{\frac{1}{2}},~~~g(x)=c_1+c_2x^{-\frac{1}{2}}$.
 
 {\bf 4)} $f(x)=x^2,~~~g(x)=c_1+c_2x$.
 
 Moreover, any admissible pair is equivalent to one of these under transformations 
 $$x\mapsto \lambda_0+\lambda_1x,~~~f\mapsto \mu_0 f+\mu_1+\mu_2x$$
  where $\lambda_1,\mu_0\ne 0$. In particular, the rank of any admissible pair is equal to 2.
 
 {\bf Proof.} In each of the cases {\bf 1), 2), 3), 4)} the corresponding cone $C$ is defined by a non-degenerate quadric. Therefore, all these cases are equivalent 
 with respect to the projective action. But the pair {\bf 4)} is admissible because it corresponds to a Gaussian integral. On the other hand, it follows from Lemma 1 that any admissible 
 pair is equivalent to one listed above. The corresponding functions $g(x)$ can be found by solving the second equation of the system (\ref{sys1}) for a known function $f(x)$. $\square$

 \section{Toward a classification of admissible pairs of functions in two variables}

 The following results are based on huge direct computations and therefore we give only a rough scheme of the proof.

 {\bf Theorem 7.1.} Let $f(x_1,x_2)$ be admissible. Then its rank is infinite iff the corresponding surface is ruled. If the corresponding surface is not ruled, then $\dim V_f\leq 6$. Moreover, $\dim V_f$ cannot be equal to 5. 

 {\bf Theorem 7.2.} Let $f(x_1,x_2)$ be admissible and $\dim V_f=6$. Then the corresponding surface is either Kummer (including degenerations of Kummer surfaces such as Steiner Roman surface) or a toric surface given by 
 $$x_0x_1x_2=x_3^3.$$

 {\bf Proof.}  Consider first three equations (\ref{f}) for $f,g$ (see Appendix) as a system of linear equations for $g$ and denote its space of solutions by $\widetilde{V}_f$. We have $V_f\subset \widetilde{V}_f$. Bringing this system of three equations to an involutive form as a system for $g$ one can see that it has an infinite-dimensional space of solutions iff $f$ corresponds to a ruled surface. Moreover, if the surface is not ruled, then $\dim \widetilde{V}_f\leq 6$ and $\dim \widetilde{V}_f=6$ only in the cases listed above. On the other hand, we know that for ruled surfaces $\dim V_f$ is infinite and for Kummer surfaces as well as for a toric variety defined above the dimension of $V_f$ is equal to 6. One can check that $\dim V_f\ne 5$ in a similar way.
 $\square$

 {\bf Remark 7.1.} It would be interesting to obtain the full classification of admissible pairs of functions in two variables. This means 

 1. Find detailed classification of admissible pairs in the case of ruled surfaces.

 2. Find all admissible pairs such that $\dim V_f=1,2,3,4$.

 {\bf Remark 7.2.} Equations (\ref{f}) become much simpler if one writes them in terms of local projective invariants. This can be done in the case $n=2$ using so-called asymptotic coordinates on a hypersurface in $\P^3$ \cite{prg}.
 It would be a good idea to use this approach for the classification of admissible pairs in two variables. 
  \section{A potential application to generalized Dirichlet series}

  In the special case when $f,g$ are both homogeneous and the equation (\ref{eq0}) can be lifted to an actual identity between distributions,  one can try to imitate the proof of the functional equation for the Riemann $\zeta$-function based on the Poisson summation formula and the Mellin transform. 
  
  Let 
  $$\int_{\R^n}g(\vec{x})e^{-\frac{f(\vec{x})}{\hbar}+\frac{i}{\hbar}\vec{x}\vec{y}}d\vec{x}=(2\pi\hbar)^{\frac{n}{2}}\hat{g}(\vec{y})e^{-\frac{\hat{f}(\vec{y})}{\hbar}}.$$
  If the Poisson summation formula is applicable (possibly after some regularization), we have
  $$\sum_{\vec{x}\in\Z^n}g(\vec{x})e^{-\frac{f(\vec{x})}{\hbar}}=(2\pi\hbar)^{\frac{n}{2}}\sum_{\vec{y}\in\Z^n}\hat{g}(2\pi\hbar\vec{y})e^{-\frac{\hat{f}(2\pi\hbar\vec{y})}{\hbar}}.$$
  Multiplying this equation by $\int_0^{\infty}\hbar^{s-1}d\hbar$, integrating by $\hbar$ and making the change of variable $t=\frac{1}{\hbar}$ we get in the l.h.s.
  $$\sum_{\vec{x}\in\Z^n}\int_0^{\infty}\hbar^{s-1} g(\vec{x})e^{-\frac{f(\vec{x})}{\hbar}}d\hbar=\sum_{\vec{x}\in\Z^n}\int_0^{\infty}t^{-1-s}g(\vec{x})e^{-f(\vec{x})t}dt=\Gamma(-s)\sum_{\vec{x}\in\Z^n}g(\vec{x})f(\vec{x})^s,$$
  and in the l.h.s. similarly, assuming that $\hat{f},\hat{g}$ are homogeneous, $\deg \hat{f}=a,~\deg\hat{g}=b$, and with the change of variable $t=\hbar^{a-1}$ we get
  $$(2\pi)^{\frac{n}{2}}\sum_{\vec{y}\in\Z^n}\int_0^{\infty}\hbar^{\frac{n}{2}}\hat{g}(2\pi\hbar\vec{y})e^{-\frac{\hat{f}(2\pi\hbar\vec{y})}{\hbar}}d\hbar=(2\pi)^{\frac{n}{2}+b}\sum_{\vec{y}\in\Z^n}\int_0^{\infty}\hbar^{\frac{n}{2}+b+s}\hat{g}(\vec{y})e^{-\hat{f}(2\pi\vec{y})\hbar^{a-1}}\frac{d\hbar}{\hbar}=$$
  $$\frac{(2\pi)^{\frac{n}{2}+b}}{a-1}\sum_{\vec{y}\in\Z^n}\hat{g}(\vec{y})\int_0^{\infty}t^{\frac{\frac{n}{2}+b+s}{a-1}}e^{-\hat{f}(2\pi\vec{y})t}\frac{dt}{t}=\frac{(2\pi)^{\frac{n}{2}+b}}{a-1}\Gamma\Bigg(\frac{\frac{n}{2}+b+s}{a-1}\Bigg)\sum_{\vec{y}\in\Z^n}\hat{g}(\vec{y})\hat{f}(2\pi\vec{y})^{-\frac{\frac{n}{2}+b+s}{a-1}}.$$
  Equating the results of these calculations for the l.h.s. and the r.h.s. we get
  $$\Gamma(-s)\sum_{\vec{x}\in\Z^n}g(\vec{x})f(\vec{x})^s=\frac{(2\pi)^{\frac{n}{2}+b}}{a-1}\Gamma\Bigg(\frac{\frac{n}{2}+b+s}{a-1}\Bigg)\sum_{\vec{y}\in\Z^n}\hat{g}(\vec{y})\hat{f}(2\pi\vec{y})^{-\frac{\frac{n}{2}+b+s}{a-1}}.$$

  The above arguments do not make rigorous sense for non-smooth functions, and one has to find a way to regularize sums and integrals above and get potentially an example of a functional equation possibly with some correction terms related with singularities. 

  For example, it would be interesting to extract a functional equation from the actual identity in Schwartz space $S^{\prime}(\R^2)$ of distributions of moderate growth:

  Let $\hbar>0$ and 
  $$\phi(x_1,x_2)=\frac{1}{(x_1-i0)^{\frac{3}{2}}}e^{\frac{i}{\hbar}(x_1+\sqrt{+ix_1}\sqrt{-ix_2})^2}.$$ 
  This function is the boundary value of a holomorphic function in $x_1$ with $\text{Im}~ x_1<0$ and in $x_2$ with $\text{Im}~ x_2>0$. Then
  $$\iint_{\R^2}\phi(x_1,x_2)e^{\frac{i}{\hbar}(x_1y_1+x_2y_2)}dx_1dx_2=2\pi i\hbar~\overline{\phi(y_2,y_1+y_2)}.$$

  \section{Conjectures and open questions}

  Let $\Sigma\subset\P^{n+1}$ be a projective hypersurface (may be a non-algebraic germ), non-degenerate at the generic point, and $C\subset\A^{n+2}$ the corresponding cone. Assume that $n\geq 2$.

  {\bf 1.} It would be interesting to determine when $rk(C)$ is infinite in terms of projective differential geometry of $\Sigma$. We know the answer in the simplest nontrivial case $n=2$: $rk(C)$ is infinite iff $\Sigma$ is a ruled surface, see Section 5.3. Notice that if $\Sigma\subset\P^3$ is a ruled surface, then its projective dual $\widehat{\Sigma}\subset\P^3$ is also ruled but the property of being ruled is not self dual if $n>2$. For arbitrary $n\geq 2$ we can only suggest the following

  {\bf Conjecture 9.1.} If $rk(C)$ is infinite, then both $\Sigma\subset\P^{n+1}$ and $\widehat{\Sigma}\subset\P^{n+1}$ are ruled. 

  Notice that the Segre cubic has rank 24 and ruled, but its dual is not ruled, see Section 5.9.

  {\bf 2.} Let $rk(C)$ be finite. It would be interesting to understand which values $rk(C)$ can take. For example, if $n=2$, then we know examples with $rk(C)=4,~6$, we also know that $rk(C)\ne 5$ and $rk(C)\leq 6$.

  {\bf Conjecture 9.2.} If $rk(C)$ is finite, then $rk(C)\leq (n+1)!$.

  It would be interesting to classify all $C$ with the finite largest possible $rk(C)$ for given $n$.

   {\bf Conjecture 9.3.} If $rk(C)$ is finite, then $\Sigma$ is algebraic.

  {\bf 3.} Let $\Sigma$ be algebraic and $rk(C)$ finite. 

  {\bf Conjecture 9.4.} If $rk(C)>0$ and finite, then $\Sigma$ is singular.

  In all {\it interesting} known examples $\Sigma$ has only isolated singularities (double points). The only exception is described in Section 5.7, these families  of surfaces have one-dimensional singularity locus. 

  Let ${\mathfrak S}_{n,d}$ be the set of algebraic hypersurfaces $\Sigma\subset\P^{n+1}$ of degree $d>2$ with the largest possible number of double points. For some values of $n,d$ all hypersurfaces from ${\mathfrak S}_{n,d}$ are admissible (with rank $(n+1)!$). For example ${\mathfrak S}_{2,4}$ are Kummer surfaces (16 double points), and ${\mathfrak S}_{3,3}$ is the Segre cubic (10 double points). It would be interesting to determine for which values of $n,d$ elements of ${\mathfrak S}_{n,d}$ are admissible.

  {\bf 4.} Let $\Sigma\subset\P^{n+1}$ be algebraic and defined by irreducible polynomial $F_1(x_0,...,x_{n+1})$, and its projective dual $\widehat{\Sigma}$ defined by an irreducible polynomial $F_2(y_0,...,y_{n+1})$. Recall that the natural birational isomorphism $\sigma:~\Sigma\to \widehat{\Sigma}$ is given by $y_i=\partial_{x_i}F_1,~i=0,...,n+1$ (see also Section 4,3, Remark 4.3.1). Let  $H=\det ((\partial_{x_i}\partial_{x_j}F_1)_{0\leq i,j\leq n+1})$. Notice that $H$ is the determinant of the Jacobian of $\sigma$, and therefore the divisor $\cal D$ on $\Sigma$ defined by $H=0$ is the singularity locus of $\sigma$. In certain examples of admissible hypersurfaces $\Sigma$ the corresponding $D$-module is holonomic, has regular singularities on $\cal D$, and $\pi_1(\Sigma\setminus{\cal D})$ acts on the space of admissible $g$. See for example Section 5.9. It would be interesting to study this class of admissible hypersurfaces. See also Section 4.5.

  {\bf 5.} The family of surfaces of degree four in Section 5.7, Theorem 5.7.2 was obtained as all possible deformations of the toric surface 
   $$x_0x_1x_2^2=x_3^4$$
   in class of surfaces of rank 4. It would be interesting to study deformations of other admissible toric hypersurfaces preserving its rank. We have checked that the toric surface 
   $$x_0x_1x_3=x_4^3$$
   with rank 6 (see Section 5.8) does not have such deformations.

  {\bf 6.} It would be interesting to study admissible pairs of the form
  $$f(x_1,...,x_N)=\sum_{i=0}^{N-m}\phi(x_{i+1},...,x_{i+m}),~~~g(x_1,...,x_N)=\sum_{i=0}^{N-m}\psi(x_{i+1},...,x_{i+m}),~~~N\gg m$$
  and similar with variables $x_{i_1,...,i_k}$, $k>1$, $1\leq i_1,...,i_k\leq N$. Such admissible pairs could be regarded as "integrable lattice models" of certain type. 

  {\bf 7.} It would be interesting to prove or disprove the following

  {\bf Conjecture 9.5.} Let $f(x_1,...,x_n)$ be a function in $n$ variables such that its Hessian is not identically zero and such that 
  $$\int g(\vec{x},\hbar)e^{\frac{f(\vec{x})+\vec{x}\cdot\vec{y}}{\hbar}}d\vec{x}=(2\pi \hbar)^{\frac{n}{2}} \hat{g}(\vec{y},\hbar)e^{\frac{\hat{f}(\vec{y})}{\hbar}}$$
  where $g(\vec{x},\hbar),~\hat{g}(\vec{y},\hbar)$ are both non-zero polynomial in $\hbar$. Then function $f$ is admissible.

  This conjecture is based on some calculation in the case $n=1$ and is not supported by any calculations for $n\geq 2$. Notice that if $g$ does not depend on $\hbar$ and $\hat{g}$ is a polynomial in $\hbar$, then pair $f,g$ is not necessarily admissible even in the case $n=1$.

 \addcontentsline{toc}{section}{Appendix. Explicit formulas for equations}

 \section*{Appendix. Explicit formulas for equations}

 Recall that $f,g$ are functions in variables $\vec{x}=(x_1,x_2,\dots,x_n)$ and we assume that the Hessian matrix $\partial^2 f:=(\partial_i \partial_j f)_{1\leq i,j \leq n}$ is non-degenerate. Denote by 
$$ (p^{ij})_{1\leq i,j\leq n}:=(\partial^2 f)^{-1}$$
the inverse matrix-valued function.

\

Main notation: for $k\ge 1$ (all summation variables and indices are assumed to be  integers),

\begin{multline*}A_k:=\sum_{v\in[0,2k]} \sum_{\substack{    d_0\in[0,\infty);\\d_1,\dots,d_v\in [3,\infty)\\
\mbox{\small such that }\\
\sum_i d_i=2(k+v),\\
d_1\geq d_2\geq \dots \geq d_v}} \sum_{\substack{a_{ij}\in [0,\infty)\\
\mbox{\small where }{0\leq i\leq j\leq v},\\
\mbox{\small satisfying }\forall i:\\
d_i=\sum_{j<i} a_{ji}+\\
+2a_{ii}+\sum_{j>i}a_{ij}}} \frac{(-1)^v}{Sym_{(d_i),(a_{ij})}}\sum_{\substack{ b_{ijl},c_{ijl}\in[1,n]\\
\mbox{\small where }\\
1\leq i\leq j\leq v\\
\mbox{\small and } 1\leq l\leq a_{ij} }}  \left(\prod_{\substack{i,j\in[1,v]\\
l\in[1, a_{ij}] \\\mbox{\small such that}\\
i\leq j}}p^{b_{ijl},c_{ijl}}\right)\cdot\\
\cdot\left[\prod_{l_1\in[1, a_{00}]} (\partial_{b_{00l_1}}\partial_{c_{00l_1}})\prod_{\substack{i\in[1,v]\\ l_2 \in [1,a_{0i}]
}}\partial_{b_{0il_2}}g\right]\cdot\prod_{i\in [1,v]}\left[\prod_{\substack{j\in [0,i)\\  l_1\in [1,a_{ji}]}}\partial_{c_{jil_1}} \prod_{l_2\in[1,a_{ii}]}(\partial_{b_{iil_2}}\partial_{c_{iil_2}})\prod_{\substack{j\in(i ,v]\\ l_3\in [1,a_{ij}]}} \partial_{b_{ijl_3}}f\right]\end{multline*}

Here  the symmetry factor is defined by
$$Sym_{(d_i),(a_{ij})}=\prod_i m_i!\cdot \prod_{\substack{i,j\in[0,v]\\
\mbox{\small such that}\\
i\leq j}}{a_{ij}!}\cdot \prod_{i\in[0,v]}{2^{a_{ii}}}$$
where  $m_1,m_2,\dots\geq 1$ are multiplicities of the repeating terms in sequence $(d_1,d_2,\dots,d_v)$, i.e.
$$d_1=\cdots=d_{m_1}>d_{m_1+1}=\cdots =d_{m_1+m_2}>d_{m_1+m_2+1}=\cdots $$

\

Meaning: let us expand $f$ at some point $\vec{x}^{(0)}$ as
$$ f=f_0+f_1+f_2+f_{\geq 3} $$
where $f_0=f(\vec{x}^{(0)})$ is a constant, $f_1,f_2$ are homogeneous polynomials in $\vec{x}-\vec{x}^{(0)}$ of degree 1 and 2 respectively, and $f_{\geq 3}$ is a series in $\vec{x}-\vec{x}^{(0)}$ containing terms of degrees $\geq 3$ only.

Then the formal Fourier transform, at point 
$$\vec{y}^{(0)}=\partial f_{|\vec{x}^{(0)}}:=(\partial_1 f,\dots,\partial_n f)_{|\vec{x}^{(0)}}$$ is equal, after normalization, to 
\begin{multline*} \int g \,e^{-\frac{f-f_0-f_1}{\hbar}} d^n\vec{x}=\int g \,e^{-\frac{f_2+f_{\geq 3}}{\hbar}} d^n\vec{x}=\\=
\sum_{v\geq 0} \frac{(-1)^v}{v!}\hbar^{-v}\int g f_{\geq 3}^v\,e^{-\frac{f_2}{\hbar}} d^n\vec{x}=\\
=(2\pi \hbar)^{n/2}\det(\partial^2 f_{|\vec{x}^{(0)}})^{-1/2}\left(g_{|\vec{x}^{(0)}}+\sum_{k\geq 1} \hbar^k {A_k}_{|\vec{x}^{(0)}}\right)\end{multline*}

In terms of (not connected) Feynman graphs, $v\ge 0$ denotes the number of vertices at which we put Taylor coefficients of $f_{\geq 3}$ (and at exactly one exceptional vertex we put Taylor coefficients of $g$). We label vertices by $\{0,1,\dots,v\}=[0,v]\cap\Z$ where $0$ corresponds to $g$, and the rest to $f_{\geq 3}$.  Moreover, we assume that the ordering of vertices is chosen in such a way that $d_1\geq d_2\geq \cdots\geq d_v\geq 3$ where for all $i\in [0,v]$ number $d_i$ is the  degree (valency) of vertex labeled by $i$. Denote by $a_{ij}\geq 0$ the number of edges connecting vertices $i$ and $j$. We enumerate edges connecting $i$ and $j$ by $\{1,\dots, a_{ij}\}$. Then we put two space indices $b_{ijl},c_{ijl}\in[1,n]$ on two ends of the edge corresponding to $l\in [1,a_{ij}]$. The factors $\prod_i m_i!$,  $\prod_{ij}a_{ij}!$ and $\prod_i 2^{a_{ii}}$ come from  symmetry, the rest is the usual Wick formula.

The total number of edges $e$ satisfies constraints:
$$  e\ge {3\over 2}v, \quad k=e-v\implies e\in \{k,\dots,3k\}\,,$$
hence in the expression $A_k$ the propagator $(p^{ij})_{1\le i,j\le n}$ appears at most $3k$ times. 

\

One-loop exactness is equivalent to an infinite sequence of differential equations:
\begin{equation}\label{f}
 A_1=0,~A_2=0,\dots
 \end{equation}

Up to symmetry, the number of distinct  graphs for $A_1,A_2,A_3$ is $5, 41, 378$ respectively.

For example, 5 graphs appearing in $A_1$ are the following:

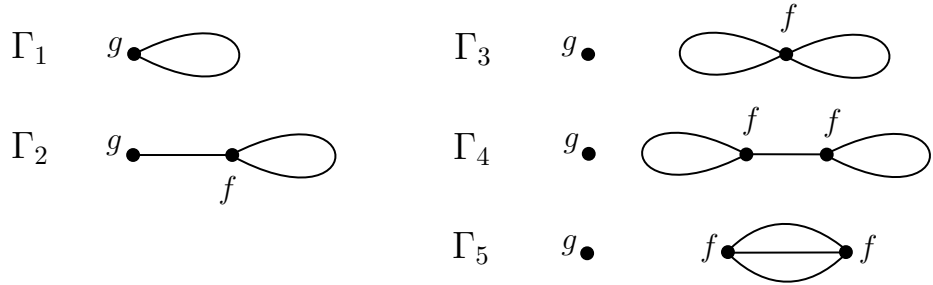
\tikzset{every picture/.style={line width=0.75pt}} 

\begin{tikzpicture}[x=0.75pt,y=0.75pt,yscale=-1,xscale=1]

\draw    (90,110.5) -- (140,110.5) ;
\draw    (139,110.6) .. controls (210.35,73.59) and (209.15,149.19) .. (139,110.5) ;
\draw  [fill={rgb, 255:red, 0; green, 0; blue, 0 }  ,fill opacity=1 ] (87,110.5) .. controls (87,108.84) and (88.34,107.5) .. (90,107.5) .. controls (91.66,107.5) and (93,108.84) .. (93,110.5) .. controls (93,112.16) and (91.66,113.5) .. (90,113.5) .. controls (88.34,113.5) and (87,112.16) .. (87,110.5) -- cycle ;
\draw  [fill={rgb, 255:red, 0; green, 0; blue, 0 }  ,fill opacity=1 ] (137,110.5) .. controls (137,108.84) and (138.34,107.5) .. (140,107.5) .. controls (141.66,107.5) and (143,108.84) .. (143,110.5) .. controls (143,112.16) and (141.66,113.5) .. (140,113.5) .. controls (138.34,113.5) and (137,112.16) .. (137,110.5) -- cycle ;
\draw  [fill={rgb, 255:red, 0; green, 0; blue, 0 }  ,fill opacity=1 ] (87.5,59.5) .. controls (87.5,57.84) and (88.84,56.5) .. (90.5,56.5) .. controls (92.16,56.5) and (93.5,57.84) .. (93.5,59.5) .. controls (93.5,61.16) and (92.16,62.5) .. (90.5,62.5) .. controls (88.84,62.5) and (87.5,61.16) .. (87.5,59.5) -- cycle ;
\draw  [fill={rgb, 255:red, 0; green, 0; blue, 0 }  ,fill opacity=1 ] (446.5,159.5) .. controls (446.5,157.84) and (447.84,156.5) .. (449.5,156.5) .. controls (451.16,156.5) and (452.5,157.84) .. (452.5,159.5) .. controls (452.5,161.16) and (451.16,162.5) .. (449.5,162.5) .. controls (447.84,162.5) and (446.5,161.16) .. (446.5,159.5) -- cycle ;
\draw  [fill={rgb, 255:red, 0; green, 0; blue, 0 }  ,fill opacity=1 ] (387,159.5) .. controls (387,157.84) and (388.34,156.5) .. (390,156.5) .. controls (391.66,156.5) and (393,157.84) .. (393,159.5) .. controls (393,161.16) and (391.66,162.5) .. (390,162.5) .. controls (388.34,162.5) and (387,161.16) .. (387,159.5) -- cycle ;
\draw  [fill={rgb, 255:red, 0; green, 0; blue, 0 }  ,fill opacity=1 ] (316.5,59.67) .. controls (316.5,61.23) and (317.84,62.5) .. (319.5,62.5) .. controls (321.16,62.5) and (322.5,61.23) .. (322.5,59.67) .. controls (322.5,58.11) and (321.16,56.85) .. (319.5,56.85) .. controls (317.84,56.85) and (316.5,58.11) .. (316.5,59.67) -- cycle ;
\draw    (402.5,110.1) -- (440.9,110.18) ;
\draw    (418.4,60.3) .. controls (489.75,23.29) and (488.55,98.89) .. (418.4,60.2) ;
\draw  [fill={rgb, 255:red, 0; green, 0; blue, 0 }  ,fill opacity=1 ] (316.8,109.9) .. controls (316.8,108.24) and (318.14,106.9) .. (319.8,106.9) .. controls (321.46,106.9) and (322.8,108.24) .. (322.8,109.9) .. controls (322.8,111.56) and (321.46,112.9) .. (319.8,112.9) .. controls (318.14,112.9) and (316.8,111.56) .. (316.8,109.9) -- cycle ;
\draw  [fill={rgb, 255:red, 0; green, 0; blue, 0 }  ,fill opacity=1 ] (416.4,59.6) .. controls (416.4,57.94) and (417.74,56.6) .. (419.4,56.6) .. controls (421.06,56.6) and (422.4,57.94) .. (422.4,59.6) .. controls (422.4,61.26) and (421.06,62.6) .. (419.4,62.6) .. controls (417.74,62.6) and (416.4,61.26) .. (416.4,59.6) -- cycle ;
\draw    (392,159.8) -- (447.6,159.68) ;
\draw    (418.4,59.7) .. controls (347.35,22.29) and (348.95,95.89) .. (418.4,59.6) ;
\draw  [fill={rgb, 255:red, 0; green, 0; blue, 0 }  ,fill opacity=1 ] (316,160.1) .. controls (316,158.44) and (317.34,157.1) .. (319,157.1) .. controls (320.66,157.1) and (322,158.44) .. (322,160.1) .. controls (322,161.76) and (320.66,163.1) .. (319,163.1) .. controls (317.34,163.1) and (316,161.76) .. (316,160.1) -- cycle ;
\draw  [fill={rgb, 255:red, 0; green, 0; blue, 0 }  ,fill opacity=1 ] (396.3,110.2) .. controls (396.3,108.54) and (397.64,107.2) .. (399.3,107.2) .. controls (400.96,107.2) and (402.3,108.54) .. (402.3,110.2) .. controls (402.3,111.86) and (400.96,113.2) .. (399.3,113.2) .. controls (397.64,113.2) and (396.3,111.86) .. (396.3,110.2) -- cycle ;
\draw    (438.83,110.37) .. controls (510.18,73.36) and (508.98,148.96) .. (438.83,110.27) ;
\draw  [fill={rgb, 255:red, 0; green, 0; blue, 0 }  ,fill opacity=1 ] (436.83,110.27) .. controls (436.83,108.61) and (438.18,107.27) .. (439.83,107.27) .. controls (441.49,107.27) and (442.83,108.61) .. (442.83,110.27) .. controls (442.83,111.92) and (441.49,113.27) .. (439.83,113.27) .. controls (438.18,113.27) and (436.83,111.92) .. (436.83,110.27) -- cycle ;
\draw    (399.4,110.2) .. controls (328.35,72.79) and (329.95,146.39) .. (399.4,110.1) ;
\draw    (391.5,157.5) .. controls (408.6,141.18) and (428.6,140.68) .. (450.6,159.18) ;
\draw    (391.1,161.18) .. controls (409.6,178.68) and (430.1,179.18) .. (449.1,160.68) ;
\draw    (90.5,59.6) .. controls (161.85,22.59) and (160.65,98.19) .. (90.5,59.5) ;

\draw (27,46.9) node [anchor=north west][inner sep=0.75pt]  [font=\large]  {$\Gamma _{1}$};
\draw (27,97.4) node [anchor=north west][inner sep=0.75pt]  [font=\large]  {$\Gamma _{2}$};
\draw (250.5,47.4) node [anchor=north west][inner sep=0.75pt]  [font=\large]  {$\Gamma _{3}$};
\draw (250,97.9) node [anchor=north west][inner sep=0.75pt]  [font=\large]  {$\Gamma _{4}$};
\draw (249.5,147.9) node [anchor=north west][inner sep=0.75pt]  [font=\large]  {$\Gamma _{5}$};
\draw (75,98.9) node [anchor=north west][inner sep=0.75pt]    {$g$};
\draw (75.5,49.4) node [anchor=north west][inner sep=0.75pt]    {$g$};
\draw (304.5,48.9) node [anchor=north west][inner sep=0.75pt]    {$g$};
\draw (305.5,97.9) node [anchor=north west][inner sep=0.75pt]    {$g$};
\draw (304.5,149.4) node [anchor=north west][inner sep=0.75pt]    {$g$};
\draw (131.5,119.4) node [anchor=north west][inner sep=0.75pt]    {$f$};
\draw (414.5,32.9) node [anchor=north west][inner sep=0.75pt]    {$f$};
\draw (395.5,84.9) node [anchor=north west][inner sep=0.75pt]    {$f$};
\draw (436.5,85.9) node [anchor=north west][inner sep=0.75pt]    {$f$};
\draw (374.5,148.9) node [anchor=north west][inner sep=0.75pt]    {$f$};
\draw (454.5,149.4) node [anchor=north west][inner sep=0.75pt]    {$f$};

\end{tikzpicture}

and the expression $A_1$ is

\begin{multline*}   {1\over 2}\sum_{i,j}p^{ij}\,\d_{ij} g-{1\over 2}\sum_{i_1 j_1 i_2 j_2}p^{i_1 j_1}p^{i_2 j_2}\,\d_{i_1}g\, \d_{j_1 i_2 j_2}f-\\
-{1\over 8}\,g\sum_{i_1 j_2 i_2 j_2}p^{i_1 j_1}p^{i_2 j_2}\,\d_{i_1 j_1 i_2 j_2}f+{1\over 8}\,g\sum_{i_1 j_2 i_2 j_2 i_3 j_3}p^{i_1 j_1}p^{i_2 j_2}p^{i_3 j_3}\,\d_{i_1 j_1 i_2}f \,\d_{j_2 i_3 j_3}f+\\
+{1\over 12}\,g\sum_{i_1 j_1 i_2 j_2 i_3 j_3}p^{i_1 j_1}p^{i_2 j_2}p^{i_3 j_3}\,\d_{i_1 i_2 i_3}f\,\d_{j_1 j_2 j_3}f\,.
\end{multline*}

\addcontentsline{toc}{section}{Acknowledgements}

\section*{Acknowledgements}

We are grateful to Robert Bryant and Joseph M. Landsberg for useful discussions. We are grateful to Nikolai Perkhunkov for useful advises and help with managing huge Maple computations.
A.O. is grateful to IHES for invitations and an excellent working atmosphere.

\end{document}